\documentclass{amsart}
\usepackage{amssymb}
\usepackage{epsfig}
\newtheorem{theorem}{Theorem}[section]
\newtheorem{lemma}[theorem]{Lemma}
\newtheorem{proposition}[theorem]{Proposition}
\newtheorem{corollary}[theorem]{Corollary}
\newtheorem{definition}[theorem]{Definition}
\newtheorem{conj}[theorem]{Conjecture}

\newcommand{\C}{\mathbb{C}}
\newcommand{\Z}{\mathbb{Z}}
\newcommand{\CP}{\mathbb{CP}}
\newcommand{\R}{\mathbb{R}}

\newcommand{\dbar}{\bar{\partial}}
\newcommand{\comment}[1]{}
\title{Lefschetz pencils, branched covers and symplectic invariants}
\author{Denis Auroux}
\address{
Department of Mathematics, M.I.T., Cambridge MA 02139, USA}
\email{
auroux@math.mit.edu}
\author{Ivan Smith}
\address{Centre for Mathematical Sciences, University of Cambridge,
Wilberforce Road, Cambridge CB3 0WB, U.K.}
\email{i.smith@dpmms.cam.ac.uk}

\begin{document}

\maketitle

\begin{center}\textit{
Lectures given at the CIME summer school ``Symplectic 4-manifolds and
algebraic surfaces'', Cetraro (Italy), September 2--10, 2003.
}\end{center}
\vskip8mm


\hfill
\begin{minipage}{0.8\linewidth}
\small
Two symplectic fibrations are never {\it exactly} the same.
When you have  two fibrations, they might be canonically 
isomorphic, but when you look closely, the points of one might
be numbers while the points of the other are bananas.

\hfill {\it (P. Seidel, 9/9/03)}
\end{minipage}
\bigskip

\section{Introduction and background}

This set of lectures aims to give an overview of Donaldson's theory of
linear systems on symplectic manifolds and the algebraic and geometric
invariants to which they give rise. After collecting some of the relevant
background, we discuss topological, algebraic and symplectic viewpoints
on Lefschetz pencils and branched covers of the projective plane. The later
lectures discuss invariants obtained by combining this theory with
pseudo-holomorphic curve methods.

\subsection{Symplectic manifolds}

\begin{definition}
A {\em symplectic structure} on a smooth manifold $M$ is a closed
non-degenerate $2$-form $\omega$, i.e.\ an element $\omega\in\Omega^2(M)$
such that $d\omega=0$ and $\forall v\in TM-\{0\}$, $\iota_v\omega\neq 0$.
\end{definition}

For example, $\R^{2n}$ carries a standard symplectic structure, given by
the 2-form $\omega_0=\sum dx_i\wedge dy_i$. Similarly, every orientable
surface is symplectic, taking for $\omega$ any non-vanishing volume form.

Since $\omega$ induces a non-degenerate antisymmetric bilinear pairing on
the tangent spaces to $M$, it is clear that every symplectic manifold is
even-dimensional and orientable (if $\dim M=2n$, then $\frac{1}{n!}\omega^n$
defines a volume form on $M$).

Two important features of symplectic structures that set them apart from
most other geometric structures are the existence of a large number of
symplectic automorphisms, and the absence of local geometric invariants.

The first point is illustrated by the following construction. Consider
a smooth function $H:M\to\R$ (a {\it Hamiltonian}), and define $X_H$ to be
the vector field on $M$ such that $\omega(X_H,\cdot)=dH$. Let
$\phi_t:M\to M$ be the family of diifeomorphisms generated by the flow of
$X_H$, i.e., $\phi_0=\mathrm{Id}$ and $\frac{d}{dt}\phi_t(x)=X_H(\phi_t(x))$.
Then $\phi_t$ is a {\it symplectomorphism}, i.e.\ $\phi_t^*\omega=\omega$.
Indeed, we have $\phi_0^*\omega=\omega$, and
$$\frac{d}{dt}\phi_t^*\omega=\phi_t^*(L_{X_H}\omega)=
\phi_t^*(d\iota_{X_H}\omega+\iota_{X_H}d\omega)=\phi_t^*(d(dH)+0)=0.$$
Therefore, the group of symplectomorphisms $\mathrm{Symp}(M,\omega)$ is
infinite-dimensional, and its Lie algebra contains all Hamiltonian vector
fields. An easy consequence is that $\mathrm{Symp}(M,\omega)$ acts
transitively on points of $M$. This is in contrast with the case of
Riemannian metrics, where isometry groups are much smaller.

The lack of local geometric invariants of symplectic structures is 
illustrated by two classical results of fundamental importance, which show
that the study of symplectic manifolds largely reduces to {\it topology}
(i.e., to discrete invariants): Darboux's theorem, and Moser's stability 
theorem. The first one shows that all symplectic forms are locally equivalent,
in sharp contrast to the case of Riemannian metrics where curvature provides
a local invariant, and the second one shows that exact deformations of
symplectic structures are trivial.

\begin{theorem}[Darboux]
Every point in a symplectic manifold $(M^{2n},\omega)$ admits a neighborhood
that is symplectomorphic to a neighborhood of the origin in
$(\R^{2n},\omega_0)$.
\end{theorem}

\proof
We first use local coordinates to map a neighborhood of a given
point in $M$ diffeomorphically onto a neighborhood $V$ of the origin in
$\R^{2n}$. Composing this diffeomorphism $f$ with a suitable linear
transformation of $\R^{2n}$, we can ensure that the symplectic form
$\omega_1=(f^{-1})^*\omega$ coincides with $\omega_0$ at the origin.
This implies that, restricting to a smaller neighborhood if necessary,
the closed 2-forms $\omega_t=t\,\omega_1 + (1-t)\omega_0$ are non-degenerate
over $V$ for all $t\in [0,1]$.

Using the Poincar\'e lemma, consider a family of $1$-forms $\alpha_t$ on $V$
such that $\frac{d}{dt}\omega_t=-d\alpha_t$. Subtracting a constant $1$-form
from $\alpha_t$ if necessary, we can assume that $\alpha_t$ vanishes at the
origin for all $t$. Using the non-degeneracy of $\omega_t$
we can find vector fields $X_t$ such that $\iota_{X_t}\omega_t = \alpha_t$.
Let $(\phi_t)_{t\in [0,1]}$ be the flow generated by $X_t$, i.e.\ the family
of diffeomorphisms defined by $\phi_0=\mathrm{Id}$, $\frac{d}{dt}\phi_t(x)=
X_t(\phi_t(x))$; we may need to restrict to a smaller neighborhood
$V'\subset V$ of the origin in order to make the flow $\phi_t$ well-defined
for all $t$. We then have
$$\frac{d}{dt}\phi_t^*\omega_t=\phi_t^*(L_{X_t}\omega_t) +
\phi_t^*\bigl(\frac{d\omega_t}{dt}\bigr)=
\phi_t^*(d(\iota_{X_t}\omega_t)-d\alpha_t)=0,$$
and therefore $\phi_1^*\omega_1=\omega_0$. Therefore, $\phi_1^{-1}\circ f$
induces a symplectomorphism from a neighborhood of $x$ in $(M,\omega)$ to a
neighborhood of the origin in $(\R^{2n},\omega_0)$.
\endproof

\begin{theorem}[Moser]
Let $(\omega_t)_{t\in [0,1]}$ be a continuous family of symplectic forms
on a compact manifold $M$. Assume that the cohomology class $[\omega_t]\in
H^2(M,\R)$ does not depend on $t$. Then $(M,\omega_0)$ is symplectomorphic
to $(M,\omega_1)$.
\end{theorem}

\proof
We use the same argument as above: since $[\omega_t]$ is constant there exist
\hbox{1-forms}
$\alpha_t$ such that $\frac{d}{dt}\omega_t=-d\alpha_t$. Define vector fields
$X_t$ such that $\iota_{X_t}\omega_t=\alpha_t$ and the corresponding flow
$\phi_t$. By the same calculation as above, we
conclude that $\phi_1^*\omega_1 = \omega_0$.
\endproof

\begin{definition} A submanifold $W\subset (M^{2n},\omega)$ is called
{\em symplectic} if $\omega_{|W}$ is non-degenerate at every point of $W$
(it is then a symplectic form on $W$);
{\em isotropic} if $\omega_{|W}=0$; and {\em Lagrangian} if it is isotropic
of maximal dimension $\dim W=n=\frac{1}{2}\dim M$.
\end{definition}

An important example is the following: given any smooth manifold $N$, the
cotangent bundle $T^*N$ admits a canonical symplectic structure that can be
expressed locally as $\omega=\sum dp_i\wedge dq_i$ (where $(q_i)$ are local
coordinates on $N$ and $(p_i)$ are the dual coordinates on the cotangent
spaces). Then the zero section is a Lagrangian submanifold of $T^*N$.

Since the symplectic form induces a non-degenerate pairing between tangent
and normal spaces to a Lagrangian submanifold, the normal bundle to a
Lagrangian submanifold is always isomorphic to its cotangent bundle. The
fact that this isomorphism extends beyond the infinitesimal level is a
classical result of Weinstein:

\begin{theorem}[Weinstein]
For any Lagrangian submanifold $L\subset (M^{2n},\omega)$, there exists
a neighborhood of $L$ which is symplectomorphic to a neighborhood of the
zero section in the cotangent bundle $(T^*L,\sum dp_i\wedge dq_i)$.
\end{theorem}

There is also a neighborhood theorem for symplectic submanifolds; in that
case, the local model for a neighborhood of the submanifold $W\subset M$
is a neighborhood of the zero section in the symplectic vector bundle
$NW$ over $W$ (since $Sp(2n)$ retracts onto $U(n)$, the classification of
symplectic vector bundles is the same as that of complex vector bundles).

\subsection{Almost-complex structures}

\begin{definition}
An {\em almost-complex structure} on a manifold $M$ is an endomorphism $J$ of
the tangent bundle $TM$ such that $J^2=-\mathrm{Id}$. An almost-complex
structure $J$ is said to be {\em tamed} by a symplectic form $\omega$ if
for every non-zero tangent vector $u$ we have $\omega(u,Ju)>0$; it is
{\em compatible} with $\omega$ if it is $\omega$-tame and 
$\omega(u,Jv)=-\omega(Ju,v)$; equivalently, $J$ is $\omega$-compatible if
and only if $g(u,v)=\omega(u,Jv)$ is a Riemannian metric.
\end{definition}

\begin{proposition}
Every symplectic manifold $(M,\omega)$ admits a compatible almost-complex
structure. Moreover, the space of $\omega$-compatible (resp.\ $\omega$-tame)
almost-complex structures is contractible.
\end{proposition}

This result follows from the fact that the space of compatible (or tame)
complex structures on a symplectic vector space is non-empty and
contractible (this can be seen by constructing explicit retractions);
it is then enough to observe that a compatible
(resp.\ tame) almost-complex structure on a symplectic manifold is simply
a section of the bundle $End(TM)$ that defines a compatible (resp.\ tame)
complex structure on each tangent space.

An almost-complex structure induces a splitting of the complexified tangent
and cotangent bundles: $TM\otimes\C=TM^{1,0}\oplus TM^{0,1}$, where
$TM^{1,0}$ and $TM^{0,1}$ are respectively the $+i$ and $-i$ eigenspaces
of $J$ (i.e., $TM^{1,0}=\{v-iJv,\ v\in TM\}$, and similarly for $TM^{0,1}$;
for example, on $\C^n$ equipped with its standard complex structure, the
$(1,0)$ tangent space is generated by $\partial/\partial z_i$ and the
$(0,1)$ tangent space by $\partial/\partial \bar{z}_i$.
Similarly, $J$ induces a complex structure on the cotangent bundle, and
$T^*M\otimes\C=T^*M^{1,0}\oplus T^*M^{0,1}$ (by definition $(1,0)$-forms
are those which pair trivially with $(0,1)$-vectors, and vice versa).
This splitting of the cotangent bundle induces a splitting of differential
forms into ``types'': $\bigwedge^rT^*M\otimes \C=
\bigoplus_{p+q=r}\bigwedge^pT^*M^{1,0} \otimes \bigwedge^qT^*M^{0,1}$.
Moreover, given a function $f:M\to\C$ we can write $df=\partial f + \dbar
f$, where $\partial f=\frac{1}{2}(df-i\,df\circ J)$ and $\dbar f=\frac{1}{2}
(df+i\,df\circ J)$ are the $(1,0)$ and $(0,1)$ parts of $df$ respectively.
Similarly, given a complex vector bundle $E$ over $M$ equipped with a
connection, the covariant derivative $\nabla$ can be split into operators
$\partial^\nabla:\Gamma(E)\to\Omega^{1,0}(E)$ and
$\dbar{}^\nabla:\Gamma(E)\to\Omega^{0,1}(E)$.

Although the tangent space to a symplectic manifold $(M,\omega)$ equipped
with a compatible almost-complex structure $J$ can be pointwise identified
with $(\C^n,\omega_0,i)$, there is an important difference between a
symplectic manifold equipped with a compatible almost-complex structure and
a complex K\"ahler manifold: the possible lack of {\it integrability} of
the almost-complex structure, namely the fact that the Lie bracket of two
$(1,0)$ vector fields is not necessarily of type $(1,0)$.

\begin{definition}
The {\em Nijenhuis tensor} of an almost-complex manifold $(M,J)$ is the
quantity defined by $N_J(X,Y)=\frac{1}{4}([X,Y]+J[X,JY]+J[JX,Y]-[JX,JY])$.
The almost-complex structure $J$ is said to be {\em integrable} if $N_J=0$.
\end{definition}

It can be checked that $N_J$ is a tensor (i.e., only depends on the values
of the vector fields $X$ and $Y$ at a given point), and that
$N_J(X,Y)=2\,\mathrm{Re}([X^{1,0},Y^{1,0}]^{(0,1)})$. The non-vanishing
of $N_J$ is therefore an obstruction to the integrability of a local frame
of $(1,0)$ tangent vectors, i.e.\ to the existence of local holomorphic
coordinates. The Nijenhuis tensor is also related to the fact that the
exterior differential of a $(1,0)$-form may have a non-zero component of
type $(0,2)$, so that the $\partial$ and $\dbar$ operators on differential
forms do not have square zero ($\dbar{}^2$ can be expressed in terms of
$\partial$ and the Nijenhuis tensor).

\begin{theorem}[Newlander-Nirenberg] Given an almost-complex manifold
$(M,J)$, the following properties are
equivalent: $(i)$~$N_J=0$; $(ii)$ $[T^{1,0}M,T^{1,0}M]\subset T^{1,0}M$;
$(iii)$ $\dbar{}^2=0$; $(iv)$ $(M,J)$ is a complex manifold, i.e.\ admits
complex analytic coordinate charts.
\end{theorem}

\subsection{Pseudo-holomorphic curves and Gromov-Witten invariants}
\leavevmode\medskip

\noindent
Pseudo-holomorphic curves, first introduced by Gromov in 1985 \cite{Gr},
have since then become the most important tool in modern symplectic
topology. In the same way as the study of complex curves in complex
manifolds plays a central role in algebraic geometry, the study of
pseudo-holomorphic curves has revolutionized our
understanding of symplectic manifolds.

The equation for holomorphic maps between two almost-complex manifolds
becomes overdetermined as soon as the complex dimension of the domain
exceeds $1$, so we cannot expect the presence of any almost-complex
submanifolds of complex dimension $\ge 2$ in a symplectic manifold equipped with
a generic almost-complex structure. On the other hand, $J$-holomorphic
curves, i.e.\ maps from a Riemann surface $(\Sigma,j)$ to the manifold
$(M,J)$ such that $J\circ df = df\circ j$ (or in the usual notation,
$\dbar_J f=0$), are governed by an elliptic
PDE, and their study makes sense even in non-K\"ahler symplectic manifolds.
The questions that we would like to answer are of the following type:
\medskip

{\it Given a compact symplectic manifold $(M,\omega)$ equipped with a generic
compatible almost-complex structure $J$ and a homology class $\beta\in H_2(M,\Z)$,
what is the number of pseudo-holomorphic curves of given genus $g$,
representing the homology class $\beta$ and passing through $r$ given points
in $M$ (or through $r$ given submanifolds in $M$)~?}%
\medskip

The answer to this question is given by {\it Gromov-Witten invariants},
which count such curves (in a sense that is not always obvious, as the
result can e.g.\ be negative, and need not even be an integer). We will
only consider a simple instance of the theory, in which we count holomorphic
spheres which are sections of a fibration.

To start with, one must study deformations of
pseudo-holomorphic curves, by linearizing the equation $\dbar_J f=0$ near a
solution. The linearized Cauchy-Riemann operator $D_{\dbar}$, whose kernel
describes infinitesimal deformations of a given curve $f:S^2\to M$,
is a Fredholm operator of (real) index
$$2d:=\mathrm{ind}\,D_{\dbar}=(\dim_\R M-6)+2\,c_1(TM)\cdot[f(S^2)].$$

When the considered curve is {\it regular}, i.e.\ when the linearized operator
$D_{\dbar}$ is surjective, the deformation theory is unobstructed,
and we expect the
moduli space $\mathcal{M}(\beta)=\{f:S^2\to M,\ \dbar_J f=0,\
[f(S^2)]=\beta\}$ to be locally a smooth manifold of
real dimension $2d$.

The main result underlying the theory of pseudo-holomorphic curves is
Gromov's compactness theorem (see \cite{Gr}, \cite{McS2}, \dots):

\begin{theorem}[Gromov]\label{thm:gromov}
Let $f_n:(\Sigma_n,j_n)\to (M,\omega,J)$ be a sequence of pseudo-holomorphic
curves in a compact symplectic manifold, representing a fixed homology
class. Then a subsequence of $\{f_n\}$
converges (in the ``Gromov-Hausdorff topology'') to a limiting map
$f_\infty$, possibly singular.
\end{theorem}

The limiting curve $f_\infty$ can have a very complicated structure, and in
particular its domain may be a nodal Riemann surface with more than one
component, due to the phenomenon of {\it bubbling}. For example, the
sequence of degree $2$ holomorphic curves $f_n:\CP^1\to\CP^2$ defined by
$f_n(u\!:\!v)=(u^2\!:\!uv\!:\!\frac{1}{n}v^2)$ converges to a singular
curve with two degree 1 components: for $(u\!:\!v)\neq
(0\!:\!1)$, we have $\lim f_n(u\!:\!v)=(u\!:\!v\!:\!0)$, so that the
sequence apparently converges to a line in $\CP^2$. However the derivatives
of $f_n$ become unbounded near $(0\!:\!1)$, and composing $f_n$ with the
coordinate change $\phi_n(u\!:\!v)=(\frac{1}{n}u\!:\!v)$ we obtain
$f_n\circ \phi_n(u\!:\!v)=(\frac{1}{n^2}u^2\!:\!\frac{1}{n}uv\!:\!
\frac{1}{n}v^2)=(\frac{1}{n}u^2\!:\!uv\!:\!v^2)$, which
converges to $(0\!:\!u\!:\!v)$ everywhere except at
$(1\!:\!0)$, giving the other component (the ``bubble'') in the
limiting curve. Therefore, it can happen that the moduli space
$\mathcal{M}(\beta)$ is not compact, and needs to be compactified by
adding maps with singular (nodal, possibly reducible) domains.

In the simplest case where the dimension of the moduli space is $2d=0$,
and assuming regularity, we can obtain an invariant by counting the
number of points of the compactified moduli space 
$\overline{\mathcal{M}}(\beta)$ (up to sign).
In the situations we consider, the moduli space will always be smooth and
compact, but may have the wrong (excess) dimension, consisting of curves
whose deformation theory is obstructed. In this case there is an {\it
obstruction bundle} $\mathrm{Obs}\to \mathcal{M}(\beta)$, whose fiber at
$(f:S^2\to M)$ is $\mathrm{Coker}\,D_{\dbar}$. In this case the invariant
may be recovered as the Euler class of this bundle, viewed as an integer
(the degree of $\mathrm{Obs}$).

\subsection{Lagrangian Floer homology}

Roughly speaking, Floer homology is a refinement of intersection theory for
Lagrangian submanifolds, in which we can only cancel intersection points by
Whitney moves along pseudo-holomorphic Whitney discs. Formally the
construction proceeds as follows.

Consider two compact orientable
(relatively spin) Lagrangian submanifolds $L_0$ and $L_1$ in
a symplectic manifold $(M,\omega)$ equipped with a compatible almost-complex
structure $J$. Lagrangian Floer homology corresponds to the Morse theory of
a functional on (a covering of) the space of arcs joining $L_0$ to $L_1$,
whose critical points are constant paths. 

For simplicity, we will only consider situations where it is not necessary
to keep track of relative homology classes (e.g.\ by working over a Novikov
ring), and where no bubbling can occur. For example, if we assume that
$\pi_2(M)=\pi_2(M,L_i)=0$, then Floer homology is well-defined; to get
well-defined product structures we will only work with exact Lagrangian
submanifolds of exact symplectic manifolds (see the final section).

By definition, the {\it Floer complex} $CF^*(L_0,L_1)$ is
the free module with one generator for each intersection
point $p\in L_0\cap L_1$, and grading given by the
{\it Maslov index}. 

Given two points $p_{\pm}\in L_0\cap L_1$, we can define a moduli space
$\mathcal{M}(p_-,p_+)$ of pseudo-holomorphic maps $f:\R\times [0,1]\to M$ 
such that $f(\cdot,0)\in L_0$, $f(\cdot,1)\in L_1$, and $\lim_{t\to \pm\infty}
f(t,\tau)=p_{\pm}$ $\forall \tau\in [0,1]$;
the expected dimension of this moduli space is the difference of Maslov
indices. Assuming regularity and compactness of $\mathcal{M}(p_-,p_+)$,
we can define an operator $\partial$ on $CF^*(L_0,L_1)$ by the formula
$$\partial p_-=\sum_{p_+}
\#(\mathcal{M}(p_-,p_+)/\R)\ p_+,$$
where the sum runs over all $p_+$ for which the expected dimension of
the moduli space is $1$.

In good cases we have $\partial^2=0$, which allows us to define
the Floer homology 
$HF^*(L_0,L_1)=\mathrm{Ker}\,\partial/\mathrm{Im}\,\partial$.
The assumptions made above on $\pi_2(M)$ and $\pi_2(M,L_i)$ eliminate the
serious technical difficulties associated to bubbling (which are more
serious than in the compact case, since bubbling can also occur on the
boundary of the domain, a real codimension 1 phenomenon which may prevent
the compactified moduli space from carrying a fundamental class, see
\cite{FO3}).

When Floer homology is well-defined, it has important consequences on the
intersection properties of Lagrangian submanifolds. Indeed, for every
{\it Hamiltonian} diffeomorphism $\phi$ we have 
$HF^*(L_0,L_1)=HF^*(L_0,\phi(L_1))$; and if $L_0$ and $L_1$ intersect
transversely, then the total rank of $HF^*(L_0,L_1)$ gives a lower bound on
the number of intersection points of $L_0$ and $L_1$. A classical
consequence, using the definition of Floer homology and the relation between
$HF^*(L,L)$ and the usual cohomology $H^*(L)$, is the non-existence of
compact simply connected Lagrangian submanifolds in $\mathbb{C}^n$.

Besides a differential, Floer complexes for Lagrangians are also equipped
with a product structure, i.e.\ a morphism $CF^*(L_0,L_1)\otimes
CF^*(L_1,L_2)\to CF^*(L_0,L_2)$ (well-defined in the cases that we will
consider). This product structure is defined as follows:
consider three points $p_1\in L_0\cap L_1$, $p_2\in L_1\cap L_2$,
$p_3\in L_0\cap L_2$, and the
moduli space $\mathcal{M}(p_1,p_2,p_3)$ of all pseudo-holomorphic
maps $f$ from a disc with three marked points $q_1,q_2,q_3$ on its
boundary to $M$, taking $q_i$ to $p_i$ and the three portions of boundary
delimited by the marked points to $L_0,L_1,L_2$ respectively.
We compactify this moduli
space and complete it if necessary in order to obtain a well-defined
fundamental cycle. The virtual
dimension of this moduli space is the difference between the Maslov index
of $p_3$ and the sum of
those of $p_1$ and $p_2$ . The product of 
$p_1$ and $p_2$ is then defined as
\begin{equation}\label{eq:hfproduct}
p_1\cdot p_2=\sum_{p_3} \#\mathcal{M}(p_1,p_2,p_3)\,
p_3,\end{equation}
where the sum runs over all $p_3$ for which the expected dimension
of the moduli space is zero.

While the product structure on $CF^*$ defined by (\ref{eq:hfproduct})
satisfies the Leibniz
rule with respect to the differential $\partial$ (and hence descends to a
product structure on Floer homology), it differs from usual products
by the fact that it is only associative {\it up to homotopy}. In fact, Floer
complexes come equipped with a full set of {\it higher-order products}
$$\mu^n:CF^*(L_0,L_1)\otimes \cdots\otimes CF^*(L_{n-1},L_n)\to CF^*(L_0,L_n)
\quad \mbox{for all}\ n\ge 1,$$
with each $\mu^n$ shifting degree by $2-n$. The first two maps $\mu^1$ and
$\mu^2$ are respectively the Floer differential $\partial$ and the
product described above. The definition of $\mu^n$ is similar to
those of $\partial$ and of the product structure: given
generators $p_i\in CF^*(L_{i-1},L_i)$ for $1\le i\le n$ and
$p_{n+1}\in CF^*(L_0,L_n)$ such that $\deg p_{n+1}=
\sum_{i=1}^n \deg p_i+2-n$, the coefficient of $p_{n+1}$
in $\mu^n(p_1,\dots,p_n)$ is obtained by counting (in a
suitable sense) pseudo-holomorphic 
maps $f$ from a disc with $n+1$ marked points $q_1,\dots,q_{n+1}$ on its
boundary to $M$, such that $f(q_i)=p_i$ and the portions of boundary
delimited by the marked points are mapped to $L_0,\dots,L_n$ respectively.

The maps $(\mu^n)_{n\ge 1}$ define an {\em $A_\infty$-structure} on Floer
complexes, i.e.\ they satisfy an infinite sequence of algebraic relations:
\smallskip
\begin{equation*}\begin{cases}
\mu^1(\mu^1(a))=0,\\
\mu^1(\mu^2(a,b))=\mu^2(\mu^1(a),b)+(-1)^{\deg a}\mu^2(a,\mu^1(b)),\\
\mu^1(\mu^3(a,b,c))=\mu^2(\mu^2(a,b),c)-\mu^2(a,\mu^2(b,c))\\
\qquad\qquad\qquad\qquad
\pm\mu^3(\mu^1(a),b,c)\pm\mu^3(a,\mu^1(b),c)\pm\mu^3(a,b,\mu^1(c)),\\
\cdots
\end{cases}\end{equation*}

This leads to the concept of ``Fukaya category'' of a symplectic manifold.
Conjecturally, for every symplectic manifold $(M,\omega)$ one should be able
to define an \hbox{$A_\infty$-category} $\mathcal{F}(M)$ whose objects are
Lagrangian submanifolds (compact, orientable, relatively spin, ``twisted''
by a flat unitary vector bundle);
the space of morphisms between two objects $L_0$ and $L_1$ is the Floer
complex $CF^*(L_0,L_1)$ equipped with its differential $\partial=\mu^1$, with
(non-associative) composition given by the product $\mu^2$, and higher order
compositions $\mu^n$.

The importance of Fukaya categories in modern symplectic topology is largely
due to the {\it homological mirror symmetry} conjecture, formulated by
Kontsevich. Very roughly, this conjecture states that the phenomenon of mirror
symmetry, i.e.\ a conjectural correspondence between symplectic manifolds
and complex manifolds (``mirror pairs'') arising from a duality among string
theories, should be visible at the level of
Fukaya categories of symplectic manifolds and categories of coherent sheaves
on complex manifolds: given a mirror pair consisting of a symplectic
manifold $M$ and a complex manifold $X$, the derived categories
$D\mathcal{F}(M)$ and $D^b Coh(X)$ should be equivalent (in a more precise form
of the conjecture, one should actually consider families of manifolds and
deformations of categories). However, due to the very
incomplete nature of our understanding of Fukaya categories in comparison to
the much better understood derived categories of coherent sheaves, this
conjecture has so far only been verified on very specific examples.

\subsection{The topology of symplectic 4-manifolds}

To end our introduction, we mention some of the known results and open
questions in the theory of compact symplectic 4-manifolds, which motivate
the directions taken in the later lectures.

Recall that, in the case of open manifolds,
Gromov's $h$-principle implies that the existence of an almost-complex
structure is sufficient. In contrast, the case of compact manifolds is much
less understood, except in dimension 4.
Whereas the existence of a class $\alpha\in H^2(M,\R)$ such that
$\alpha^{\cup n}\neq 0$ and of an almost-complex structure already provide
elementary obstructions to the existence of a symplectic structure on a
given compact manifold, in the case of 4-manifolds a much stronger
obstruction arises from Seiberg-Witten invariants. We will not define
these, but mention some of their key topological consequences for symplectic
4-manifolds, which follow from the work of Taubes (\cite{Ta1}, \cite{Ta2},
\dots).

\begin{theorem}[Taubes]\label{thm:taubes}
$(i)$ Let $(M^4,\omega)$ be a compact symplectic 4-manifold with $b_2^+\ge 2$.
Then the homology class $c_1(K_M)$ admits a (possibly disconnected) smooth
pseudo-holo\-mor\-phic representative (in particular $c_1(K_M)\cdot
[\omega]\ge 0$).
Hence, if $M$ is minimal i.e.\ contains no $(-1)$-spheres,
then $c_1(K_M)^2=2\chi(M)+3\sigma(M)\ge 0$.

$(ii)$ If $(M^4,\omega)$ splits as a connected sum $M_1\#M_2$, then one
of the $M_i$ has negative definite intersection form.
\end{theorem}

When $b_2^+(M)=1$, Seiberg-Witten theory still has some implications.
Using Gromov's characterization of the Fubini-Study symplectic structure
of $\CP^2$ in terms of the existence of
pseudo-holomorphic lines, Taubes has shown that the symplectic structure of
$\CP^2$ is unique up to scaling. This result has been extended by Lalonde
and McDuff to the case of rational ruled surfaces, where $\omega$ is
determined by its cohomology class.
\medskip

{\bf Remark.} 
For any smooth connected symplectic curve
$\Sigma$ in a symplectic four-manifold $(M,\omega)$, the genus
$g(\Sigma)$ is related to the homology class by the classical
{\it adjunction formula}
$$2-2g(\Sigma)+[\Sigma]\cdot[\Sigma]=-c_1(K_M)\cdot[\Sigma],$$ a direct
consequence of the splitting $TM_{|\Sigma}=T\Sigma\oplus N\Sigma$.
For example, every connected component of the pseudo-holomorphic
representative of $c_1(K_M)$ constructed by Taubes satisfies
$g(\Sigma)=1+[\Sigma]\cdot[\Sigma]$ (this is how one derives the
inequality $c_1(K_M)^2\ge 0$ under the minimality assumption).
In fact, Seiberg-Witten theory also implies that symplectic curves have
minimal genus among all smoothly embedded surfaces in their homology class.
\medskip

In parallel to the above constraints on symplectic 4-manifolds, surgery
techniques have led to many interesting new examples of compact symplectic
manifolds.

One of the most efficient techniques in this respect is the {\it symplectic
sum} construction, investigated by Gompf \cite{Go1}:
if two symplectic manifolds $(M_1^{2n},\omega_1)$ and $(M_2^{2n},\omega_2)$
contain compact symplectic hypersurfaces $W_1^{2n-2},W_2^{2n-2}$ that are
mutually symplectomorphic and whose normal bundles have opposite Chern
classes, then we can cut $M_1$ and $M_2$ open along the submanifolds $W_1$
and $W_2$, and glue them to each other along their common boundary,
performing a fiberwise connected sum in the normal bundles to $W_1$ and
$W_2$, to obtain a new symplectic manifold $M=M_1{\vphantom{M}}_{W_1}{\#}
{\vphantom{M}}_{W_2} M_2$. This construction has in particular allowed Gompf
to show that every finitely presented group can be realized as the
fundamental group of a compact symplectic 4-manifold. This is in sharp
contrast to the K\"ahler case, where Hodge theory shows that the first
Betti number is always even.

A large number of questions remain open, even concerning the Chern numbers
realized by symplectic 4-manifolds. For instance it is unknown to this date
whether the Bogomolov-Miyaoka-Yau inequality $c_1^2\le 3 c_2$, satisfied by
all complex surfaces of general type, also holds in the symplectic case.
Moreover, very little is known about the symplectic topology of 
complex surfaces of general type.

\section{Symplectic Lefschetz fibrations}

This section will provide a theoretical classification of symplectic
4-manifolds in algebraic terms, but we begin very humbly.

\subsection{Fibrations and monodromy}

Here is an easy way to build symplectic 4-manifolds:

\begin{proposition}[Thurston]
If $\Sigma_g\to X\to \Sigma_h$ is a surface bundle with fiber non-torsion
in homology, then $X$ is symplectic.
\end{proposition}

\proof
Let $\eta\in\Omega^2(M)$ be a closed 2-form representing a cohomology
class which pairs non-trivially with the fiber. Cover the base $\Sigma_h$ by
balls $U_i$ over which the fibration
is trivial: we have a diffeomorphism $\phi_i:f^{-1}(U_i)\to
U_i\times \Sigma_g$, which determines a projection $p_i:f^{-1}(U_i)\to
\Sigma_g$.

Let $\sigma$ be a symplectic form on the fiber $\Sigma_g$, in the same
cohomology class as the restriction of $\eta$.
After restriction to $f^{-1}(U_i)\simeq U_i\times \Sigma_g$, we
can write $p_i^*\sigma=\eta+d\alpha_i$ for some 1-form $\alpha_i$
over $f^{-1}(U_i)$.
Let $\{\rho_i\}$ be a partition of unity subordinate to
the cover $\{U_i\}$ of $\Sigma_h$, and let
$\tilde\eta=\eta+\sum_i d((\rho_i\circ f)\,\alpha_i)$.
The 2-form $\tilde\eta$ is well-defined since the support of $\rho_i\circ f$
is contained in $f^{-1}(U_i)$, and it is obviously closed. Moreover, over 
$f^{-1}(p)$, we have
$\tilde\eta_{|f^{-1}(p)}=\eta_{|f^{-1}(p)}+\sum_i \rho_i(p)
d\alpha_{i|f^{-1}(p)}=\sum_i \rho_i(p)\,(\eta+d\alpha_i)_{|f^{-1}(p)}=
\sum_i \rho_i(p)\,(p_i^*\sigma)_{|f^{-1}(p)}$. Since a positive linear
combination of symplectic forms over a Riemann surface is still symplectic,
the form $\tilde\eta$ is non-degenerate on every fiber.

At any point $x\in X$,
the tangent space $T_xX$ splits into a vertical subspace
$V_x=\mathrm{Ker}\,df_x$ and a horizontal subspace $H_x=\{v\in T_xX,\ 
\tilde\eta(v,v')=0\ \forall v'\in V_x\}$. Since the restriction of
$\tilde\eta$ to the vertical subspace is non-degenerate, we have
$T_xX=H_x\oplus V_x$. Letting $\kappa$ be a symplectic form on the
base $\Sigma_h$, the 2-form $f^*\kappa$ is non-degenerate over
$H_x$, and therefore for sufficiently large $C>0$ the 2-form
$\tilde\eta+C\,f^*\kappa$ defines a global symplectic form on $X$.
\endproof

The cohomology class of the symplectic form depends on $C$, but the
structure is canonical up to deformation equivalence.
The hypothesis on the fiber is satisfied whenever $g\neq 1$, since
$c_1(TX^{\mathrm{vert}})$ evaluates non-trivially on the fiber. That
some assumption is needed for $g=1$ is shown by the example of the
Hopf fibration $T^2\to S^1\times S^3\to S^2$.
Historically, the first example of a non-K\"ahler symplectic 4-manifold,
due to Thurston \cite{Th}, is a non-trivial $T^2$-bundle over $T^2$
(the product of $S^1$ with the mapping torus of a Dehn twist, which
has $b_1=3$).

Unfortunately, not many four-manifolds are fibered.

\begin{definition}
A Lefschetz pencil on a smooth oriented four-manifold $X$ is a map
$f:X-\{b_1,\dots,b_n\}\to S^2$, submersive away from a finite set
$\{p_1,\dots,p_r\}$, conforming to local models $(i)$ $(z_1,z_2)\mapsto
z_1/z_2$ near each $b_i$, $(ii)$ $(z_1,z_2)\mapsto z_1^2+z_2^2$ near each
$p_j$. Here the $z_i$ are orientation-preserving local complex-valued
coordinates.
\end{definition}

We can additionally require that the critical
values of $f$ are all distinct (so that each fiber contains at most one
singular point).

This definition is motivated by the complex analogue of Morse
theory. Global holomorphic functions on a projective surface must be
constant, but interesting functions exist on the complement of finitely
many points, and the generic such will have only quadratic singularities.
The (closures of the) fibers of the map $f$ cut the four-manifold $X$ into
a family of real surfaces all passing through the $b_i$ (locally like
complex lines through a point in $\C^2$), and with certain fibers having
nodal singularities ($(z_1+iz_2)(z_1-iz_2)=0$).
If we blow up the $b_i$, then the map $f$ extends to the entire manifold
and we obtain a {\it Lefschetz fibration}.

A small generalization of the previous argument to the case of Lefschetz
fibrations shows that, if $X$ admits
a Lefschetz pencil, then it is symplectic (work on the blow-up and 
choose the constant $C$ so large that the exceptional sections arising
from the $b_i$ are all symplectic and can be symplectically blown down).
In fact the symplectic form obtained in this way is canonical up to isotopy
rather than just deformation equivalence, as shown by Gompf \cite{GS,Go2}.

A real Morse function encodes the topology of a manifold: the critical
values disconnect $\R$, and the topology of the level sets changes by a
handle addition as we cross a critical value. In the complex case, the
critical values do not disconnect, but the local model is determined by
its {\it monodromy}, i.e.\ the diffeomorphism of the smooth fiber obtained
by restricting the fibration to a circle enclosing a single critical value.

The fiber $F_t$ of the map
$(z_1,z_2)\to z_1^2+z_2^2$ above $t\in\C$ is given by the equation
$(z_1+iz_2)(z_1-iz_2)=t$: the fiber $F_t$ is smooth 
(topologically an annulus) for all $t\neq 0$, while the fiber
above the origin presents a
transverse double point, and is obtained from the nearby fibers by
collapsing an embedded simple closed loop called the {\it vanishing cycle}.
For example, for $t>0$ the vanishing cycle is the loop
$\{(x_1,x_2)\in\R^2,\ x_1^2+x_2^2=t\}=F_t\cap \R^2\subset
F_t$.

\begin{proposition}
For a circle in the base $S^2$ of a Lefschetz fibration enclosing a single
critical value, whose
critical fiber has a single node, the monodromy is a Dehn twist about the
vanishing cycle.
\end{proposition}

\proof[Sketch of proof]
By introducing a cutoff function $\psi$ and by identifying the fiber
$z_1^2+z_2^2=t$ with the set $z_1^2+z_2^2=\psi(\|z\|^2)\,t$, we can see that
the monodromy is the identity outside a small neighborhood of the vanishing
cycle. This reduces the problem to the local model of the annulus, which
has mapping class group (relative to the boundary) isomorphic to $\Z$, so
we just need to find one integer. One possibility is to study an example,
e.g.\ an elliptic surface, where we can determine the monodromy by
considering its action on homology, interpreted as periods. Alternatively,
we can think of the annulus as a double cover of the disc branched at two
points (the two square roots of $t$), and watch these move as $t$ follows
the unit circle.

\begin{center}
\epsfig{file=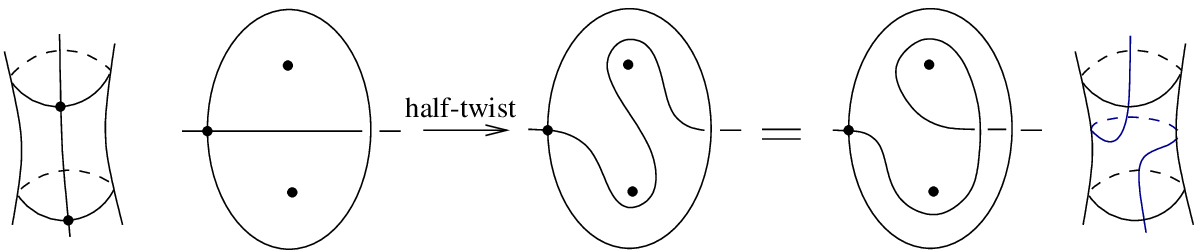}
\end{center}

One can also consider higher-dimensional Lefschetz fibrations, which by
definition are again submersions away from non-degenerate quadratic
singularities. In the local model, the smooth fibers of
$f:(z_i)\mapsto \sum z_i^2$ also contain
a (Lagrangian) sphere, and the monodromy around the critical fiber
$f^{-1}(0)$ is a generalized Dehn twist about this vanishing cycle
(see the lectures by Seidel in this volume).

Using this local model, we may now define the monodromy homomorphism.
Fix a base point $q_0\in S^2-\mathrm{crit}(f)$, and consider a closed loop
$\gamma:[0,1]\to S^2-\mathrm{crit}(f)$ (starting and ending at $q_0$). 
By fixing a horizontal distribution we can perform parallel transport in
the fibers of $f$ along $\gamma$, which induces a diffeomorphism from
$\Sigma=f^{-1}(q_0)$ to itself. (Such a horizontal distribution is
canonically defined if we fix a symplectic form on the total space by
taking the symplectic orthogonal complements to the fibers.)
The isotopy class of this diffeomorphism,
which is well-defined independently of the chosen horizontal
distribution, is called the {\it monodromy} of $f$ along $\gamma$. Hence, we
obtain a monodromy homomorphism characteristic of the Lefschetz fibration $f$,
$$\psi:\pi_1(S^2-\mathrm{crit}(f),q_0)\to
\pi_0\mathrm{Diff}^+(\Sigma),$$
which takes a loop encircling one critical value to a Dehn twist as above.
\medskip

{\bf Example:} let $C_0$ and $C_1$ be generic cubic curves in $\CP^2$, and
consider the pencil of cubics $\{C_0+\lambda C_1=0\}_{\lambda\in\CP^1}$.
This pencil has 9 base points (the intersections of $C_0$ and $C_1$), and
12 singular fibers. To see the latter fact, note that the Euler
characteristic $\chi(X)$ of the total space of a Lefschetz pencil of genus
$g$ is given by $\chi(X)=4-4g-\# b_i+\# p_j$.

After blowing up the base points, we obtain an elliptic Lefschetz fibration,
whose monodromy takes values in the genus 1 mapping class group
$\mathrm{Map}_1=\pi_0\mathrm{Diff}^+(T^2)=SL(2,\Z)$. Each local monodromy
is conjugate to $\left(\!\!\begin{array}{cc}1&1\\0&1\end{array}\!\!\right)$.
The monodromy homomorphism is determined by the images of a basis for 
$\pi_1(S^2-\mathrm{crit}(f))$ consisting of 12 loops encircling one critical
value each. The monodromy around the product of these loops is the identity,
as the product loop bounds a disc in $S^2$ over which the fibration is
trivial. In an appropriate basis, the resulting word in 12 Dehn twists in
$SL(2,\Z)$ can be brought into standard form
$$(A\cdot B)^6=\left(\left(\!\!\begin{array}{cc}1&1\\0&1\end{array}\!\!\right)\cdot
\left(\!\!\!\begin{array}{rc}1&0\\-1&1\end{array}\!\!\right)\right)^6=I.$$
Such a word in Dehn twists, called a {\it positive relation} in the relevant
mapping class group, captures the topology of a Lefschetz fibration. In the
case where the fibration admits distinguished sections, e.g.\ coming from
the base points of a pencil, we can refine the monodromy by working in the
relative mapping class group of the pair $(\Sigma,\{b_i\})$ (see the example
in the last lecture).

In fact, a careful analysis of positive relations in $SL(2,\Z)$ implies that
all elliptic Lefschetz fibrations are K\"ahler (and have monodromy words
of the form $(A\cdot B)^{6n}=1$), a classical result of Moishezon and Livne
\cite{Mo1}. More geometrically, this result can also be deduced from the
work of Siebert and Tian \cite{ST} described in their lectures in this volume.
\medskip

{\bf Remark.}
In the case of Lefschetz fibrations over a disc, the monodromy homomorphism
determines the total space of the fibration up to symplectic deformation.
When considering fibrations over $S^2$, the monodromy data determines the
fibration over a large disc $D$ containing all critical values, after
which we only need to add a trivial fibration over a small disc $D'=S^2-D$,
to be glued in a manner compatible with the fibration structure over the
common boundary $\partial D=\partial D'=S^1$. This gluing involves the
choice of a map from $S^1$ to $\mathrm{Diff}^+(\Sigma_g)$, i.e.\ an element
of $\pi_1\mathrm{Diff}^+(\Sigma_g)$, which is trivial if $g\ge 2$ (hence
in this case the positive relation determines the topology completely).
\medskip

Combining the facts that Lefschetz pencils carry symplectic structures and
that they correspond to positive relations has algebraic consequences for
the mapping class group.

\begin{proposition}
There is no positive relation involving only Dehn twists about separating
curves.
\end{proposition}

\proof[Sketch of proof] (see \cite{SmHodge} for a harder proof) ---
Suppose for contradiction we have such a word of length $\delta$.
This defines a four-manifold
$X$ with a Lefschetz fibration having this as monodromy. We can compute
the signature $\sigma(X)$ by surgery, cutting the
manifold open along neighborhoods of the singular fibers; we find (cf.\
\cite{Oz}) that each
local model contributes $-1$ to the signature, so we obtain 
$\sigma(X)=-\delta$.
This allows us to compute the Betti and Chern numbers of $X$: $b_1(X)=2g$,
$b_2(X)=\delta+2$, $c_2(X)=4-4g+\delta$, and $c_1^2(X)=8-8g-\delta<-\delta$.
Hence $c_1^2$ of any minimal model of $X$ is negative, and so $X$ must be
rational or ruled by a theorem of Liu. These cases can be excluded by hand.
\endproof

The pencil of cubics on $\CP^2$ is an instance of a much more general
construction.

\begin{proposition}[Lefschetz]
Projective surfaces have Lefschetz pencils.
\end{proposition}

Generic hyperplane sections cut out smooth complex curves, and a pencil
corresponds to a line of such hyperplanes. Inside the dual projective space
$(\CP^N)^*$ pick a line transverse to the dual variety (the locus of
hyperplanes which are not transverse to the given surface). A local
computation shows that this transversality condition goes over to give
exactly the non-degenerate critical points of a Lefschetz pencil.
 From another point of view, if $L$ is a very ample line bundle on $X$
(so sections generate the fibers of $L$), we can consider the evaluation
map $X\times H^0(L)\to L$. Let $Z\subset X\times H^0(L)$ be the preimage
of the zero section, then a regular value of the projection $Z\to H^0(L)$
is a section with smooth zero set. In this way the construction of embedded
complex curves in $X$ (and more generally linear systems of such) can be
reduced to the existence of regular values, i.e.\ Sard's theorem.

Certainly the converse to the above cannot be true: not all Lefschetz
fibrations are K\"ahler. The easiest way to see this is to use the (twisted)
fiber sum construction. Given a positive relation $\tau_1\dots\tau_s=1$ in
$\mathrm{Map}_g$, and some element $\phi\in \mathrm{Map}_g$, we obtain
a new positive relation $\tau_1\dots\tau_s(\phi^{-1}\tau_1\phi)\dots
(\phi^{-1}\tau_s\phi)=1$. If $\phi=\mathrm{Id}$ the corresponding
four-manifold is the double branched cover of the original manifold over
a union of two smooth fibers, but in general the operation has no
holomorphic interpretation. The vanishing cycles of the new fibration
are the union of the old vanishing cycles and their images under $\phi$.
Since $H_1(X)=H_1(\Sigma)/\langle\mathrm{vanishing\ cycles}\rangle$, we
can easily construct examples with odd first Betti number, for example
starting with a genus 2 pencil on $T^2\times S^2$ \cite{OS}.

More sophisticated examples (for instance with trivial first Betti number)
can be obtained by forming infinite families of twisted fiber sums with
non-conjugate monodromy groups and invoking the following

\begin{theorem}[Arakelov-Par\v{s}in]
Only finitely many isotopy classes of Lefschetz fibrations with given fiber
genus and number of critical fibers can be K\"ahler.
\end{theorem}

{\bf Remark.}
Twisted fiber sum constructions can often be ``untwisted'' by subsequently
fiber summing with another suitable (e.g.\ holomorphic) Lefschetz fibration.
A consequence of this is a stable isotopy result for genus 2 Lefschetz
fibrations \cite{Au6}: any genus 2 Lefschetz fibration becomes isotopic
to a holomorphic fibration after repeated fiber sums with the standard
holomorphic fibration with 20 singular fibers coming from a genus 2 pencil
on $\CP^1\times\CP^1$. More generally a similar result holds for all
Lefschetz fibrations with monodromy contained in the {\it hyperelliptic}
subgroup of the mapping class group. This is a corollary of a recent result of
Kharlamov and Kulikov \cite{KK} about braid monodromy factorizations: 
after repeated (untwisted) fiber sums with copies of a same fixed
holomorphic 
fibration with $8g+4$ singular fibers, any hyperelliptic genus $g$ Lefschetz
fibration eventually becomes holomorphic. Moreover, the fibration obtained
in this manner is completely
determined by its number of singular fibers of each type (irreducible,
reducible with components of given genus), and when the fiber genus is odd
by a certain $\Z_2$-valued invariant. 
(The proof of this result uses the fact that the hyperelliptic mapping
class group is an extension by $\Z_2$ of the braid group of $2g+2$ points
on a sphere, which is itself a quotient of $B_{2g+2}$; this makes it
possible to transform the monodromy of a hyperelliptic Lefschetz fibration
into a factorization in $B_{2g+2}$, with different types of factors
for the various types of singular fibers and extra contributions belonging to
the kernel of the morphism $B_{2g+2}\to B_{2g+2}(S^2)$, and hence reduce the
problem to that studied by Kharlamov and Kulikov. This connection between
mapping class groups and braid groups will be further studied in later
lectures.) It is not clear whether the result should be expected to remain
true in the non-hyperelliptic case.%
\medskip

These examples of Lefschetz fibrations differ somewhat in character from
those obtained in projective geometry, since the latter always admit
exceptional sections of square $-1$. However, an elementary argument in
hyperbolic geometry shows that fiber sums never have this property
\cite{Sm2}. To see that not every Lefschetz pencil arises from a holomorphic
family of surfaces appears to require some strictly deeper machinery.

To introduce this, let us note that, if we choose a metric on the total
space of a genus $g$ Lefschetz pencil or fibration,
the fibers become Riemann surfaces and this induces a
map $\phi:\CP^1\to\overline{\mathcal{M}}_g$ to the Deligne-Mumford moduli
space of stable genus $g$ curves. There is a line bundle
$\lambda\rightarrow\overline{\mathcal{M}}_g$
(the Hodge bundle), with fiber $\det H^0(K_{\Sigma})$ above $[\Sigma]$,
and an index theorem for the family of $\dbar$-operators on the
fibers shows

\begin{proposition}[\cite{SmHodge}]
$\sigma(X)=4\langle c_1(\lambda),[\phi(\CP^1)]\rangle-\delta$.
\end{proposition}

On the other hand, $c_1(\lambda)$ and the Poincar\'e duals of the components
of the divisor of nodal curves generate $H^2(\overline{\mathcal{M}}_g,\Z)$,
so the above formula -- together with the numbers of singular fibers of
different topological types -- characterizes the homology class
$[\phi(\CP^1)]$.

Clearly, holomorphic Lefschetz fibrations give rise to rational curves in
the moduli space, and these have locally positive intersection with all
divisors in which they are not contained. This gives another constraint
on which Lefschetz pencils and fibrations can be holomorphic. For example,
the divisor $\overline{\mathcal{H}}_3$ of hyperelliptic genus 3 curves
has homology class $[\overline{\mathcal{H}}_3]=9c_1(\lambda)-[\Delta_0]-
3[\Delta_1]$, where $\Delta_0$ and $\Delta_1$ are the divisors of irreducible
and reducible nodal curves respectively.

\begin{corollary}
A genus $3$ Lefschetz fibration $X$ with irreducible fibers and such that
$(i)$ $\chi(X)+1$ is not divisible by $7$, $(ii)$ $9\sigma(X)+5\chi(X)+40<0$
is not holomorphic.
\end{corollary}

\proof[Sketch of proof] For hyperelliptic fibrations of any genus, we have
$$(8g+4)\,c_1(\lambda)\cdot[\phi(S^2)]=g\,[\Delta_0]\cdot[\phi(S^2)]+
\sum_{h=1}^{[g/2]}4h(g-h)\,[\Delta_h]\cdot[\phi(S^2)].$$
To prove this, we can represent the four-manifold as a double cover of a
rational ruled surface (see the lectures of Siebert and Tian in this
volume). This gives another expression for $\sigma(X)$ which can be compared
to that above. Then assumption $(i)$ and integrality of the signature show
that the fibration is not isotopic to a hyperelliptic fibration, while
assumption $(ii)$ shows that $[\phi(S^2)]\cdot[\overline{\mathcal{H}}_3]<0$.
\endproof

It is possible to build a Lefschetz fibration (with 74 singular fibers)
admitting a $(-1)$-section and satisfying the conditions of the Corollary
\cite{Sm1}.

The right correspondence between the geometry and the algebra comes from
the work of Donaldson \cite{Do2,Do3}:

\begin{theorem}[Donaldson]\label{thm:donaldson}
Any compact symplectic 4-manifold admits 
Lefschetz pencils with symplectic fibers (if $[\omega]$ is integral,
Poincar\'e dual to $k[\omega]$ for any sufficiently large $k$).
\end{theorem}

As explained later, we even get some uniqueness -- but only asymptotically 
with the parameter $k$. Increasing the parameter $k$ makes the algebraic
monodromy descriptions more and more complicated, but in principle, as with
surgery theory for high-dimensional smooth manifolds, this gives a complete
algebraic encoding. The construction is flexible enough to impose some extra
conditions on the pencils, for instance we can assume that all the singular
fibers are irreducible. (If the four-manifold has even intersection form,
this is completely elementary.)

In order to describe Donaldson's construction of symplectic Lefschetz
pencils, we need a digression into
{\it approximately holomorphic geometry}.

\subsection{Approximately holomorphic geometry}

On an almost-complex manifold, the lack of integrability usually
prevents the existence of non-trivial holomorphic sections of vector bundles or
pseudo-holomorphic maps to other manifolds, but one can work in a similar
manner with approximately holomorphic objects.

Let $(M^{2n},\omega)$ be a compact symplectic manifold of dimension $2n$. We
will assume throughout this paragraph that $\frac{1}{2\pi}[\omega]\in
H^2(M,\Z)$; this integrality condition does not restrict the topological
type of $M$, since any symplectic form can be perturbed into another
symplectic form $\omega'$ whose cohomology class is rational (we can then
achieve integrality by multiplication by a constant factor). Morever, it is
easy to check that the submanifolds of $M$ that we will construct are not
only $\omega'$-symplectic but also $\omega$-symplectic, hence making the
general case of Theorem \ref{thm:donaldson} follow from the integral case.

Let $J$ be an almost-complex structure compatible with $\omega$, and
let $g(.,.)=\omega(.,J.)$ be the corresponding Riemannian metric.
We consider a complex line bundle $L$ over $M$ such that
$c_1(L)=\frac{1}{2\pi}[\omega]$, endowed with a Hermitian metric and
a Hermitian connection $\nabla^L$ with curvature 2-form 
$F(\nabla^L)=-i\omega$. The almost-complex structure induces a
splitting of the connection~:
$\nabla^L=\partial^L+\dbar^L$, where $\partial^L s(v)=\frac{1}{2}(\nabla^L
s(v)-i\nabla^L s(Jv))$ and $\dbar^L s(v)=\frac{1}{2}(\nabla^L
s(v)+i\nabla^L s(Jv))$.

If the almost-complex structure $J$ is integrable, i.e.\ if $M$ is
K\"ahler, then $L$ is an ample holomorphic line bundle,
and for large enough values of $k$ the holomorphic sections of 
$L^{\otimes k}$ determine an embedding of the manifold $M$ 
into a projective space (Kodaira's theorem). Generic hyperplane
sections of this projective embedding are smooth hypersurfaces in $M$,
and a pencil of hyperplanes through a generic codimension 2 linear subspace
defines a Lefschetz pencil.

When the manifold $M$ is only symplectic, the lack of integrability of $J$
prevents the existence of holomorphic sections. Nonetheless, it is possible
to find an {\it approximately holomorphic} local model: a neighborhood of a
point $x\in M$, equipped with the symplectic form $\omega$ and the
almost-complex structure $J$, can be identified with a neighborhood of the
origin in $\C^n$ equipped with the standard symplectic form $\omega_0$ and
an almost-complex structure of the form $i+O(|z|)$. In this local model, the
line bundle $L^{\otimes k}$ endowed with the connection
$\nabla=(\nabla^L)^{\otimes k}$ of curvature $-ik\omega$ can be identified
with the trivial line bundle $\underline{\C}$ endowed with the connection 
$d+\frac{k}{4}\sum (z_j\,d\bar{z}_j-\bar{z}_j\,dz_j)$. The section of 
$L^{\otimes k}$ given in this trivialization by $s_{\mathrm{ref},k,x}(z)=
\exp(-\frac{1}{4}k|z|^2)$ is then approximately holomorphic \cite{Do1}.

More precisely, a sequence of sections $s_k$ of $L^{\otimes k}$ is
said to be approximately holomorphic if, with respect to the rescaled
metrics $g_k=kg$, and after normalization of the sections to ensure that
$\|s_k\|_{C^r,g_k}\sim C$, an inequality of the form
$\|\dbar s_k\|_{C^{r-1},g_k}<C'k^{-1/2}$ holds, where $C$ and $C'$ are
constants independent of $k$. The change of metric, which dilates all
distances by a factor of $\smash{\sqrt{k}}$, is required in order to be able to
obtain uniform estimates, due to the large curvature of
the line bundle $L^{\otimes k}$. The intuitive idea is that, for large $k$,
the sections of the line bundle $L^{\otimes k}$ with curvature
$-ik\omega$ probe the geometry of $M$ at small scale
($\sim 1/\sqrt{k}$), which makes the almost-complex structure
$J$ almost integrable and allows one to achieve better
approximations of the holomorphicity condition $\dbar s=0$.

Since the above requirement is an open condition,
there is no well-defined ``space of approximately holomorphic
sections'' of $L^{\otimes k}$.
Nonetheless, the above local model gives us a large number of approximately
holomorphic sections (consider $s_{\mathrm{ref},k,x}$ for a large finite set of $x\in X$),
which can be used to embed $X$ as a symplectic
submanifold of a (high-dimensional) projective space. However, this embedding
by itself is not very useful since it is not clear that any of its hyperplane
sections can be used to define a smooth symplectic hypersurface in $X$.

Hence, in contrast with the complex case, the non-trivial part of the construction
is to find, among all the available approximately holomorphic sections,
some whose geometric behavior is as generic as possible. That this is at
all possible is a subtle observation of Donaldson, which leads to
the following result \cite{Do1}:

\begin{theorem}[Donaldson]\label{thm:skdhyp}
For $k\gg 0$, $L^{\otimes k}$ admits approximately holomorphic sections
$s_k$ whose zero sets $W_k$ are smooth symplectic hypersurfaces.
\end{theorem}

The proof of this result starts from the observation that, if the section $s_k$ vanishes
transversely and if $|\dbar s_k(x)|\ll |\partial s_k(x)|$ at every point of
$W_k=s_k^{-1}(0)$, then the submanifold $W_k$ is symplectic,
and even approximately $J$-holomorphic
(i.e.\ $J(TW_k)$ is close to $TW_k$). The crucial point is therefore to
obtain a lower bound for $\partial s_k$ at every point of $W_k$, in order
to make up for the lack of holomorphicity.

Sections $s_k$ of $L^{\otimes k}$ are said to be {\it
uniformly transverse to $0$} if there exists a constant $\eta>0$
(independent of $k$) such that the inequality $|\partial s_k(x)|_{g_k}>\eta$
holds at any point of $M$ where $|s_k(x)|<\eta$. In order to prove Theorem
\ref{thm:skdhyp}, it is sufficient to achieve this uniform estimate on the tranversality of
some approximately holomorphic sections $s_k$. The idea of the construction
of such sections consists of two main steps. The first one is an effective
local transversality result for complex-valued functions, for which
Donaldson's argument appeals to ideas of Yomdin about the complexity of real
semi-algebraic sets (see \cite{Au5} for a simplification of the argument).
The second step is a remarkable globalization process,
which makes it possible to achieve uniform transversality over larger and
larger open subsets by means of successive perturbations of the sections
$s_k$, until transversality holds over the entire manifold $M$ \cite{Do1}.

That the interplay between the two steps above is subtle can already be
gathered from the delicate statement of the local transversality result.  
For $\beta > 0$ set $t(\beta) =\beta/( \log\beta^{- 1})^d$. Here $\beta$
represents the maximum size of the allowed perturbation
$s_\mathrm{given} \mapsto s_\mathrm{given}-w s_\mathrm{ref}$ with $|w|<\beta$,
$t(\beta)$ is the amount of transversality thereby obtained, and $d=d(n)$ is a
universal constant that we will mostly ignore. Write $B^+ \subset \C^n$
for a Euclidean ball slightly larger than the unit ball $B$.

\begin{theorem}[\cite{Do1,Au2}] \label{transversality}
If $f:B^+\rightarrow\mathbb{C}^{n + 1}$ satisfies $| f | < 1$ and
$|\overline{\partial}f| < t(\beta)$ pointwise, then there is some
$w\in\mathbb{C}^{n+1}$
with $| w | <\beta$ such that $|f(z)-w| > t(\beta)$ over $B\subset B^+$.
\end{theorem}

To see the relevance of this, for an approximately holomorphic section $s_k$
of ${L}^{\otimes k}$ we consider the holomorphic 1-jet $(s_k,\partial s_k)$, a
section of $L^{\otimes k} \oplus L^{\otimes k} \otimes T^*M^{1,0}$. This
is locally a map from the complex $n$-dimensional ball to $\mathbb{C}^{n +
1}$, and by adding reference sections we can explicitly give local
perturbations of this. In this context, theorems such as the one above are
elementary if we take a polynomial function $t(\beta) = (const)\beta^q$ with
$q>1$, but such perturbations do not patch well.

By way of an example, let $f:\mathbb{R}^4\rightarrow\mathbb{R}^6$ have
bounded derivative, so $f$ takes balls of radius $\epsilon$ to balls of
comparable size.  Now $B^4(1)$ is filled by approximately $\epsilon^{-4}$
balls, and each is taken to a ball of volume approximately $\epsilon^6$. 
Hence, the total volume of the image is about $\epsilon^2$ of $B^6(1)$.  To
change the function $f\mapsto f - w$ by $\beta$ in order to miss the image
of $f$ by $t(\beta) =(const)\beta^3$, say, would be straightforward, because
taking $\epsilon=t(\beta)$, the $\epsilon$-neighborhood of $f(B^4)$ has a
volume of the order of $\epsilon^2$ and hence cannot contain any ball of
radius $\epsilon^{1/3}\sim \beta$.

However, our manifold is covered by $O(k^{2n})$ balls of fixed $g_k$-radius
in which our reference
sections $s_{\mathrm{ref},k}$ are concentrated.  Perturbing over each ball
one by one, all estimates are destroyed and it is impossible to achieve
uniform transversality.  The solution is to perturb over balls at great distance
simultaneously; nonetheless the simultaneous perturbations will
\emph{not} be entirely independent. We cover $X$ by a fixed number $D^{2n}$
of collections of balls (the number of balls in each collection, but not the
number of collections, will grow with the parameter $k$), in such a way that
any two balls in the same collection are at $g_k$-distance at least $D$ from
each other.  To obtain
uniformly transverse sections, we (1) start with some approximately
holomorphic section (e.g. the zero section $s_0 = 0$); (2) perturb by
$\beta_0$ over balls of the first collection $I_0$ to get a section $s_1$which
is $t(\beta_0)$ transverse over $\cup_{i\in I_0}B_i$; (3) perturb over balls
of the second collection $I_1$ by an amount $\beta_1 \ll t(\beta_0)/2$ to get a
section $ s_2$ which is $t(\beta_1)$-transverse over $\cup_{i\in I_0\cup
I_1}B_i$ etc.  Continuing, we need the sequence
$$\beta_0,\ \beta_1\sim t(\beta_0),\ \beta_2\sim t(\beta_1)\ldots$$ 
to be chosen so that, at each stage (for all
$N<D^{2n}$), $\exp(-D^2\beta_N) < \beta_{N+1}/2.$ 
Here, the left-hand side $e^{-D^2}\beta_N$ is the effect of the perturbation
at a ball $B_1$ in the collection $I_N$ on another ball $B_2$ of the
\emph{same} collection (which we perturb simultaneously); the right-hand
side $\beta_{N+1}/2$ is the transversality obtained at the ball $B_2$ by
virtue of its \emph{own} perturbation.
We have $N = D^{2n}$ stages, for some large $D$, so we need 
$$\exp(-N^{1/n})\beta_N < \beta_{N+1}/2$$ 
for $N\gg 0$.  This inequality fails for any polynomial function
$\beta_{N+1}=t(\beta_N) =\beta_N^q / 2$ for $q>1$: one gets
$\beta_N \sim (1/2)^{q^N}$ and $\exp(-N^{1/n}) \not < \beta_{N+1}/\beta_N$.

In other words, the estimates coming naively from Sard's theorem don't
provide a good enough local theorem to pass to a global one.  The remedy is
that our functions are approximately holomorphic and not arbitrary, and for
holomorphic functions stronger Sard-like theorems are available: the
prototype here is that the regular values of a smooth map $B(1)\subset \R^{2n}
\rightarrow \R^2$ are in general only dense, but the regular values of a
holomorphic map $B(1)\subset \C^n \rightarrow \C$ form the complement of a finite set. 
In practice, the proof of Theorem \ref{transversality} proceeds by reduction
to the case of polynomial functions, and then appeals either to real
algebraic geometry \cite{Do1} or to the classical monotonicity theorem
\cite{Au5}.

The symplectic submanifolds constructed by Donaldson present several
remarkable properties which make them closer to complex submanifolds than
to arbitrary symplectic submanifolds. For instance, they satisfy the
Lefschetz hyperplane theorem: up to half the dimension of the submanifold,
the homology and homotopy groups of $W_k$ are identical to those of $M$
\cite{Do1}. More importantly, these submanifolds are, in a sense,
asymptotically unique: for given large enough $k$, the submanifolds 
$W_k$ are, up to symplectic isotopy, independent of all the choices made
in the construction (including that of the almost-complex structure $J$)
\cite{Au1}.

It is worth mentioning that analogues of this
construction  have been obtained for contact manifolds by
Ibort, Martinez-Torres and Presas (\cite{IMP}, \dots\!\!); see also recent
work of Giroux and Mohsen \cite{Gi}.
\medskip

As an application of this theorem, we mention a symplectic packing result
due to Biran \cite{Biran}. Recall that a full symplectic packing of a manifold
$(M,\omega)$ is an embedding $(\amalg B_i,\omega_{std})\hookrightarrow
(M,\omega)$ of a disjoint union of standard Euclidean symplectic balls of
equal volumes whose images fill the entire volume of $M$. Gromov pointed
out that, in contrast to the volume-preserving case, there are obstructions
to symplectic packing: for instance, $\CP^2$ cannot be fully packed by two
balls.

\begin{theorem}[Biran]
Let $(M^4,\omega)$ be a symplectic four-manifold with integral $[\omega]$.
Then $M$ admits a full packing by $N$ balls for all large $N$.
\end{theorem}

The key ingredient in the proof is to reduce to the case of ruled surfaces
by decomposing $M$ into a disc bundle over a Donaldson submanifold $\Sigma$
dual to $k[\omega]$ and an isotropic CW-complex -- which takes no volume.
It is therefore sufficient to fully pack a ruled surface by balls which
remain disjoint from a section at infinity. On the other hand, there is a
well-known correspondence between embeddings of symplectic balls of size
$\mu$ and symplectic forms on a blow-up giving each exceptional curve area
$\mu$. Moreover, symplectic forms in a given cohomology class can be
constructed by symplectic inflation in the presence of appropriate embedded
symplectic surfaces, and for ruled surfaces these are provided by an
elementary computation of the Gromov or Seiberg-Witten invariants. The
best value of $N$ is not known in general, because it is determined by 
the best value of $k$ in Donaldson's construction.
\medskip

We now move on to Donaldson's construction of symplectic Lefschetz pencils
\cite{Do2,Do3}. In comparison with Theorem \ref{thm:skdhyp}, the general setup is the
same, the main difference being that we consider no longer one, but
two sections of $L^{\otimes k}$. A pair of suitably chosen approximately
holomorphic sections $(s_k^0,s_k^1)$ of $L^{\otimes k}$ defines a family of
symplectic hypersurfaces $$\Sigma_{k,\alpha}=\{x\in M,\ s_k^0(x)-\alpha
s_k^1(x)=0\},\ \ \alpha\in\CP^1=\C\cup\{\infty\}.$$ The submanifolds 
$\Sigma_{k,\alpha}$ are all smooth except for finitely many of them
which present an isolated singularity; they intersect transversely along
the {\it base points} of the pencil, which form a smooth symplectic
submanifold $Z_k=$\hbox{$\{s_k^0=s_k^1=0\}$} of codimension $4$.

The two sections $s_k^0$ and $s_k^1$ determine a projective map
$f_k=(s_k^0:s_k^1):M-Z_k\to\CP^1$, whose critical points correspond to the
singularities of the fibers $\Sigma_{k,\alpha}$. In the case of a symplectic
Lefschetz pencil, the function $f_k$ is a
complex Morse function, i.e.\ near any of its critical points it is given
by the local model $f_k(z)=z_1^2+\dots+z_n^2$ in approximately holomorphic
coordinates. After blowing up $M$ along $Z_k$, the Lefschetz pencil
structure on $M$ gives rise to a well-defined map $\hat{f}_k:\hat{M}\to
\CP^1$; this map is a symplectic Lefschetz fibration.
Hence, Theorem \ref{thm:donaldson} may be reformulated more precisely as follows:

\begin{theorem}[Donaldson]
For large enough $k$, the given manifold $(M^{2n},\omega)$ admits symplectic
Lefschetz pencil structures determined by pairs of suitably chosen
approximately holomorphic sections $s_k^0,s_k^1$ of $L^{\otimes k}$.
Moreover, for large enough $k$ these Lefschetz pencil structures are
uniquely determined up to isotopy.
\end{theorem}

As in the case of submanifolds, Donaldson's argument relies on successive
perturbations of given approximately holomorphic sections $s_k^0$ and
$s_k^1$ in order to achieve uniform transversality properties, not only for
the sections $(s_k^0,s_k^1)$ themselves but also for the derivative
$\partial f_k$ \cite{Do3}.

The precise meaning of the uniqueness statement is the following:
assume we are given two sequences of Lefschetz pencil structures on $(M,\omega)$,
determined by pairs of approximately holomorphic sections of $L^{\otimes k}$
satisfying uniform transversality estimates, but possibly with respect to
two different $\omega$-compatible almost-complex structures on $M$. Then,
beyond a certain (non-explicit) value of $k$, it becomes possible to find
one-parameter families of Lefschetz pencil structures interpolating between
the given ones. In particular, this implies that for large $k$ the monodromy
invariants associated to these Lefschetz pencils only depend on
$(M,\omega,k)$ and not on the choices made in the construction.

The monodromy invariants associated to a symplectic Lefschetz pencil are
essentially those of the symplectic Lefschetz fibration obtained after
blow-up along the base points, with only a small refinement. After the
blow-up operation, each fiber of $\hat{f_k}:\hat{M}\to \CP^1$ contains a
copy of the base locus $Z_k$ embedded as a smooth symplectic hypersurface.
This hypersurface lies away from all vanishing cycles, and is preserved by
the monodromy. Hence, the monodromy homomorphism can be defined to take
values in the group of isotopy classes of symplectomorphisms of the fiber
$\Sigma_k$ whose restriction to the submanifold $Z_k$ is the identity.

\section{Symplectic branched covers of $\CP^2$}

\subsection{Symplectic branched covers}

\begin{definition}\label{def:cover}
A smooth map $f:X^4\to (Y^4,\omega_Y)$ from a compact oriented smooth
4-manifold to a compact symplectic 4-manifold is a (generic)
{\em symplectic branched
covering} if, given any point $p\in X$, there exist 
neighborhoods $U\ni p$ and $V\ni f(p)$ and orientation-preserving local
diffeomorphisms $\phi:U\to\C^2$ and $\psi:V\to\C^2$, such that
$\psi_*\omega_Y(v,iv)>0$ $\forall v\neq 0$ $($i.e.\ the standard complex
structure is $\psi_*\omega_Y$-tame$)$, and such that
$\psi\circ f\circ \phi^{-1}$ is one of the following model maps:

$(i)$ $(u,v)\mapsto (u,v)$ (local diffeomorphism),

$(ii)$ $(u,v)\mapsto (u^2,v)$ (simple branching),

$(iii)$ $(u,v)\mapsto (u^3-uv,v)$ (cusp).
\end{definition}

The three local models appearing in this definition are exactly those
describing a generic holomorphic map between complex surfaces, except that
the local coordinate systems we consider are not holomorphic.

By computing the Jacobian of $f$ in the given local coordinates, we can see
that the {\it ramification curve} $R\subset X$ is a smooth submanifold (it
is given by $\{u=0\}$ in the second local model and $\{v=3u^2\}$ in the
third one). However, the image $D=f(R)\subset X$ (the {\it branch curve},
or {\it discriminant curve}) may be singular. More precisely, in the simple
branching model $D$ is given by $\{z_1=0\}$, while in the cusp model we have
$f(u,3u^2)=(-2u^3,3u^2)$, and hence $D$ is locally identified with
the singular curve $\{27z_1^2=4z_2^3\}\subset\C^2$. This means that, at the
cusp points, $D$ fails to be immersed. Besides the cusps, the branch curve
$D$ also generically presents {\it transverse double points} (or {\it nodes}),
which do not appear in the local models because they correspond to simple
branching in two distinct points $p_1,p_2$ of the same fiber of $f$. 
There is no constraint on the orientation of the local intersection between
the the two branches of $D$ at a node (positive or negative,
i.e.\ complex or anti-complex), because the local models near $p_1$ and
$p_2$ hold in different coordinate systems on $Y$.

Generically, the only singularities of the branch curve
$D\subset Y$ are transverse double points (``nodes'') of either orientation
and complex cusps. Moreover, because the local
models identify $D$ with a complex curve, the tameness condition on the
coordinate systems implies that $D$ is a (singular) symplectic submanifold
of~$Y$.

The following result states that a symplectic branched cover of a symplectic
4-manifold carries a natural symplectic structure \cite{Au2}:

\begin{proposition}\label{prop:cover}
If $f:X^4\to (Y^4,\omega_Y)$ is a symplectic branched cover, then $X$
carries a symplectic form $\omega_X$ such that $[\omega_X]=f^*[\omega_Y]$,
canonically determined up to symplectomorphism.
\end{proposition}

\proof
The 2-form $f^*\omega_Y$ is closed, but it is only non-degenerate outside of
$R$. At any point $p$ of $R$, the 2-plane $K_p=\mathrm{Ker}\,df_p\subset
T_pX$ carries a natural orientation induced by the complex orientation in
the local coordinates of Definition \ref{def:cover}. Using the local models, we can
construct an {\it exact} 2-form $\alpha$ such that, at any point $p\in R$,
the restriction of $\alpha$ to $K_p$ is non-degenerate and positive. 

More precisely, given $p\in R$ we consider a small ball centered at $p$ and 
local coordinates $(u,v)$ such that $f$ is given by one of the models of
the definition, and we set $\alpha_p=d(\chi_1(|u|)\chi_2(|v|)\, x\,dy)$,
where $x=\mathrm{Re}(u)$, $y=\mathrm{Im}(u)$, and $\chi_1$ and $\chi_2$ are
suitably chosen smooth cut-off functions. We then define $\alpha$ to be
the sum of these $\alpha_p$ when $p$ ranges over a
finite subset of $R$ for which the supports of the $\alpha_p$ cover the
entire ramification curve $R$. 
Since $f^*\omega_Y\wedge \alpha$ is positive at every point of $R$, it
is easy to check that the 2-form $\omega_X=f^*\omega_Y + \epsilon\,\alpha$
is symplectic for a small enough value of the constant $\epsilon>0$.

The fact that $\omega_X$ is canonical up to symplectomorphism follows
immediately from Moser's stability theorem and from the observation that
the space of exact perturbations
$\alpha$ such that $\alpha_{|K_p}>0$ $\forall p\in R$ is a convex subset
of $\Omega^2(X)$ and hence connected.
\endproof

Approximately holomorphic techniques make it possible to show that every
compact symplectic 4-manifold can be realized as a branched cover of
$\CP^2$. The general setup is
similar to Donaldson's construction of symplectic Lefschetz pencils: we
consider a compact symplectic manifold $(X,\omega)$, and perturbing the
symplectic structure if necessary we may assume that
$\frac{1}{2\pi}[\omega]\in H^2(X,\Z)$. Introducing an almost-complex 
structure $J$ and a line bundle $L$ with $c_1(L)=\frac{1}{2\pi}[\omega]$,
we consider triples of approximately holomorphic sections
$(s_k^0,s_k^1,s_k^2)$ of $L^{\otimes k}$: for $k\gg 0$, it is again possible
to achieve a generic behavior for the projective map
$f_k=(s_k^0:s_k^1:s_k^2):X\to\CP^2$ associated with the
linear system. If the manifold $X$ is four-dimensional, then the linear
system generically has no base points, and for a suitable choice of sections
the map $f_k$ is a branched covering \cite{Au2}.

\begin{theorem}\label{thm:cover}
For large enough $k$, three suitably chosen approximately holomorphic
sections of $L^{\otimes k}$ over $(X^4,\omega)$ determine a symplectic
branched covering $f_k:X^4\to\CP^2$, described in approximately holomorphic
local coordinates by the local models of Definition \ref{def:cover}.
Moreover, for $k\gg
0$ these branched covering structures are uniquely determined up to isotopy.
\end{theorem}

Because the local models hold in approximately holomorphic (and hence
$\omega$-tame) coordinates, the ramification curve $R_k$ of $f_k$ is a
symplectic submanifold in $X$ (connected, since the Lefschetz hyperplane
theorem applies). Moreover, if we normalize the Fubini-Study symplectic form
on $\CP^2$ in such a way that $\frac{1}{2\pi}[\omega_{FS}]$ is the generator
of $H^2(\CP^2,\Z)$, then we have $[f_k^*\omega_{FS}]=2\pi c_1(L^{\otimes k}=
k[\omega]$, and it is fairly easy to check that the symplectic form on $X$
obtained by applying Proposition \ref{prop:cover} to the branched covering $f_k$
coincides up to symplectomorphism with $k\omega$ \cite{Au2}. In fact, the
exact 2-form
$\alpha=k\omega-f_k^*\omega_{FS}$ is positive over $\mathrm{Ker}\,df_k$ at
every point of $R_k$, and $f_k^*\omega_{FS}+t\alpha$ is a symplectic form
for all $t\in (0,1]$.

The uniqueness statement in Theorem \ref{thm:cover}, which should be interpreted exactly
in the same way as that obtained by Donaldson for Lefschetz pencils, implies
that for $k\gg 0$ it is possible to define invariants of the symplectic
manifold $(X,\omega)$ in terms of the monodromy of the branched covering
$f_k$ and the topology of its branch curve $D_k\subset \CP^2$. However,
the branch curve $D_k$ is only determined up to creation or cancellation
of (admissible) pairs of nodes of opposite orientations.

A similar construction can be attempted when $\dim X>4$; in this case, the
set of base points $Z_k=\{s_k^0=s_k^1=s_k^2=0\}$ is no longer empty; it
is generically a smooth codimension 6 symplectic submanifold.
With this understood, Theorem \ref{thm:cover} admits the following higher-dimensional
analogue \cite{Au3}:

\begin{theorem}
For large enough $k$, three suitably chosen approximately holomorphic
sections of $L^{\otimes k}$ over $(X^{2n},\omega)$ determine a map
$f_k:X-Z_k\to\CP^2$ with generic local models, canonically determined up to
isotopy.
\end{theorem}

The model maps describing the local behavior of $f_k$ in approximately
holomorphic local coordinates are now the following:\smallskip

\begin{tabular}{cl}
$(0)$ & $(z_1,\dots,z_n)\mapsto (z_1:z_2:z_3)$ near a base point,\\
$(i)$ & $(z_1,\dots,z_n)\mapsto (z_1,z_2)$,\\
$(ii)$ & $(z_1,\dots,z_n)\mapsto (z_1^2+\dots+z_{n-1}^2,z_n)$,\\
$(iii)$ & $(z_1,\dots,z_n)\mapsto (z_1^3-z_1z_n+z_2^2+\dots+z_{n-1}^2,z_n)$.
\end{tabular}
\smallskip

\noindent
The set of critical points $R_k\subset X$ is again a (connected) smooth
symplectic curve, and its image $D_k=f_k(R_k)\subset\CP^2$ is again a
singular symplectic curve whose only singularities generically are
transverse double points of either orientation and complex cusps.
The fibers of $f_k$ are codimension 4 symplectic submanifolds, intersecting
along $Z_k$; the fiber above a point of $\CP^2-D_k$ is smooth, while
the fiber above a smooth point of $D_k$ presents an ordinary double point,
the fiber above a node presents two ordinary double points,
and the fiber above a cusp presents an $\mathrm{A}_2$ singularity.

The proof of these two results relies on a careful examination of the
various possible local behaviors for the map $f_k$ and on transversality
arguments establishing the existence of sections of $L^{\otimes k}$ with generic
behavior. Hence, the argument relies on the enumeration of the various
special cases, generic or not, that may occur; each one corresponds to the
vanishing of a certain quantity that can be expressed in terms of the
sections $s_k^0,s_k^1,s_k^2$ and their derivatives. Therefore, the proof
largely reduces to a core ingredient which imitates classical singularity
theory and can be thought of as a uniform transversality result for jets of
approximately holomorphic sections \cite{Au4}.

Given approximately holomorphic sections $s_k$ of very positive bundles
$E_k$ (e.g.\ $E_k=\C^m\otimes L^{\otimes k}$) over the symplectic manifold
$X$, one can consider the $r$-jets $j^r s_k=(s_k,\partial s_k,(\partial
\partial s_k)_\mathrm{sym}, \dots, (\partial^r s_k)_\mathrm{sym})$, which
are sections of the {\it jet bundles} $\mathcal{J}^r E_k=\bigoplus_{j=0}^r
(T^*X^{(1,0)})_{\mathrm{sym}}^{\otimes j}\otimes E_k$. Jet bundles can
naturally be stratified by approximately holomorphic submanifolds
corresponding to the various possible local behaviors at order $r$ for the
sections $s_k$. The generically expected behavior corresponds to the case
where the jet $j^r s_k$ is transerse to the submanifolds in the
stratification. The result is the following \cite{Au4}:

\begin{theorem}
Given stratifications $\mathcal{S}_k$ of the jet bundles $\mathcal{J}^r E_k$
by a finite number of approximately holomorphic submanifolds
(Whitney-regular, uniformly transverse to fibers, and with curvature bounded
independently of $k$), for large enough $k$ the vector bundles $E_k$ admit
approximately holomorphic sections $s_k$ whose $r$-jets are uniformly
transverse to the stratifications $\mathcal{S}_k$. Moreover these sections
may be chosen arbitrarily close to given sections.
\end{theorem}

A one-parameter version of this result also holds, which makes it
possible to obtain results of asymptotic uniqueness up to isotopy for
generic sections.

Applied to suitably chosen stratifications, this result provides the main
ingredient for the construction of $m$-tuples of approximately
holomorphic sections of $L^{\otimes k}$ (and hence projective maps $f_k$ to
$\CP^{m-1}$) with generic behavior. Once uniform transversality
of jets has been obtained, the only remaining task is to achieve some
control over the antiholomorphic derivative $\dbar f_k$ near the critical
points of $f_k$ (typically its vanishing in some directions), in order to
ensure that $\dbar f_k \ll \partial f_k$ everywhere; for low values of $m$
such as those considered above, this task is comparatively easy.

\subsection{Monodromy invariants for branched covers of $\CP^2$}

The topological data characterizing a symplectic branched covering
$f:X^4\to\CP^2$ are on one hand the topology of the branch curve
$D\subset\CP^2$ (up to isotopy and cancellation of pairs of nodes),
and on the other hand a monodromy morphism $\theta:\pi_1(\CP^2-D)\to S_N$
describing the manner in which the $N=\deg f$ sheets of the covering are
arranged above $\CP^2-D$.

Some simple properties of the monodromy morphism $\theta$ can be readily
seen by considering the local models of Definition \ref{def:cover}. For example, the
image of a small loop $\gamma$ bounding a disc that intersects $D$
transversely in a single smooth point (such a loop is called a {\it geometric
generator} of $\pi_1(\CP^2-D)$) by $\theta$ is necessarily a transposition.
The smoothness of $X$ above a singular point of $D$ implies some
compatibility properties on these transpositions (geometric generators
corresponding to the two branches of $D$ at a node must map to disjoint
commuting transpositions, while to a cusp must correspond a pair of adjacent
transpositions). Finally, the connectedness of $X$ implies the surjectivity
of $\theta$ (because the subgroup $\mathrm{Im}(\theta)$ is generated by
transpositions and acts transitively on the fiber of the covering).

It must be mentioned that the amount of information present in the monodromy
morphism $\theta$ is fairly small: a classical conjecture in algebraic
geometry (Chisini's conjecture, essentially solved by Kulikov \cite{Ku})
asserts that, given an algebraic singular plane curve $D$ with cusps and nodes, a
symmetric group-valued monodromy morphism $\theta$ compatible with $D$
(in the above sense), if it exists, is unique except for a small list of
low-degree counter-examples.
Whether Chisini's conjecture also holds for symplectic branch
curves is an open question, but in any case the number of possibilities for
$\theta$ is always finite.

The study of a singular complex curve $D\subset\CP^2$ can be carried out
using the braid monodromy techniques developed in complex algebraic geometry
by Moishezon and Teicher \cite{Mo2,Te1,...}: the idea is to choose a linear
projection $\pi:\CP^2-\{\mathrm{pt}\}\to\CP^1$, for example
$\pi(x\!:\!y\!:\!z)=(x\!:\!y)$, in such a way that the curve $D$ lies in
general position with respect to the fibers of $\pi$, i.e.\ $D$ is
positively transverse to the fibers of $\pi$ everywhere except at isolated
non-degenerate smooth complex tangencies. The restriction $\pi_{|D}$ is then
a singular branched covering of degree $d=\deg D$, with {\it special points}
corresponding to the singularities of $D$ (nodes and cusps) and to the
tangency points. Moreover, we can assume that all special points lie in
distinct fibers of $\pi$. A plane curve satisfying these topological
requirements is called a {\it braided} (or {\it Hurwitz}) curve.
\medskip

\begin{center}
\setlength{\unitlength}{0.8mm}
\begin{picture}(80,55)(-40,-15)
\put(0,-2){\vector(0,-1){8}}
\put(2,-7){$\pi:[x\!:\!y\!:\!z]\mapsto [x\!:\!y]$}
\put(-40,-15){\line(1,0){80}}
\put(-38,-12){ $\CP^1$}
\put(-40,0){\line(1,0){80}}
\put(-40,40){\line(1,0){80}}
\put(-40,0){\line(0,1){40}}
\put(40,0){\line(0,1){40}}
\put(-38,33){ $\CP^2-\{\infty\}$}
\put(27,31){ $D$}
\multiput(-20,20)(0,-2){18}{\line(0,-1){1}}
\multiput(-5,20)(0,-2){18}{\line(0,-1){1}}
\multiput(15,15)(0,-2){9}{\line(0,-1){1}}
\multiput(15,-9)(0,-2){3}{\line(0,-1){1}}
\put(-20,-15){\circle*{1}}
\put(-5,-15){\circle*{1}}
\put(15,-15){\circle*{1}}
\qbezier[200](25,35)(5,30)(-5,20)
\qbezier[60](-5,20)(-10,15)(-15,15)
\qbezier[40](-15,15)(-20,15)(-20,20)
\qbezier[40](-20,20)(-20,25)(-15,25)
\qbezier[60](-15,25)(-10,25)(-5,20)
\qbezier[110](-5,20)(0,15)(15,15)
\qbezier[320](15,15)(5,15)(-30,5)
\put(-20,20){\circle*{1}}
\put(-5,20){\circle*{1}}
\put(15,15){\circle*{1}}
\end{picture}
\end{center}
\medskip

Except for those which contain special points of $D$, the fibers of $\pi$
are lines intersecting the curve $D$ in $d$ distinct points. If one chooses
a reference point $q_0\in\CP^1$ (and the corresponding fiber $\ell\simeq\C
\subset\CP^2$ of $\pi$), and if one restricts to an affine subset in order
to be able to trivialize the fibration $\pi$, the topology of the branched
covering $\pi_{|D}$ can be described by a {\it braid monodromy} morphism
\begin{equation}
\rho:\pi_1(\C-\{\mathrm{pts}\},q_0)\to B_d,
\end{equation}
where $B_d$ is the braid group on $d$ strings. The braid $\rho(\gamma)$
corresponds to the motion of the $d$ points of $\ell\cap D$ inside the
fibers of $\pi$ when moving along the loop $\gamma$.

Recall that the braid group $B_d$ is the fundamental group of the
configuration space of $d$ distinct points in $\R^2$; it is also the
group of isotopy classes of compactly supported orientation-preserving
diffeomorphisms of $\R^2$ leaving invariant a set of $d$ given distinct
points. It is generated by the standard {\it half-twists} $X_1,\dots,X_{d-1}$
(braids which exchange two consecutive points by rotating them 
counterclockwise by 180 degrees around each other), with relations
$X_iX_j=X_jX_i$ for $|i-j|\ge 2$ and $X_iX_{i+1}X_i=X_{i+1}X_iX_{i+1}$
(the reader is referred to Birman's book \cite{Bi} for more details).

Another equivalent way to consider the monodromy of a braided curve is
to choose an ordered system of generating
loops in the free group $\pi_1(\C-\{\mathrm{pts}\},q_0)$. The morphism
$\rho$ can then be described by a {\it factorization} in the braid group
$B_d$, i.e.\ a decomposition of the monodromy at infinity into the product
of the individual monodromies around the various special points of $D$.
By observing that the total space of $\pi$ is the line bundle $O(1)$ over
$\CP^1$, it is easy to see that the monodromy at infinity is given by the
central element $\Delta^2=(X_1\dots X_{d-1})^d$ of $B_d$ (called ``full
twist'' because it represents a rotation of a large disc by 360 degrees).
The individual monodromies around the special points are conjugated to
powers of half-twists, the exponent being $1$ in the case of tangency
points, $2$ in the case of positive nodes (or $-2$ for negative nodes), and
$3$ in the case of cusps.

The braid monodromy $\rho$ and the corresponding factorization depend on
trivialization choices, which affect them by {\it simultaneous conjugation} by
an element of $B_d$ (change of trivialization of the fiber $\ell$ of $\pi$),
or by {\it Hurwitz operations} (change of generators of the group
$\pi_1(\C-\{\mathrm{pts}\},q_0)$). There is a one-to-one correspondence between
braid monodromy morphisms $\rho:\pi_1(\C-\{\mathrm{pts}\})\to B_d$ (mapping
generators to suitable powers of half-twists) up to
these two algebraic operations and singular (not necessarily complex)
braided curves of degree $d$ in $\CP^2$ 
up to isotopy among such curves (see e.g.\ \cite{KK} for a detailed
exposition). Moreover, it is easy to check that
every braided curve in $\CP^2$ can be deformed into a braided symplectic
curve, canonically up to isotopy among symplectic braided curves (this
deformation is performed by collapsing the curve $D$ into a neighborhood
of a complex line in a way that preserves the fibers of $\pi$). However,
the curve $D$ is isotopic to a complex curve only for certain specific
choices of the morphism $\rho$.

Unlike the case of complex curves, it is not clear {\it a priori} that the
symplectic branch curve $D_k$ of one of the covering maps given by
Theorem \ref{thm:cover} can be made compatible with the linear projection $\pi$;
making the curve $D_k$ braided relies on an improvement of the result in
order to control more precisely the behavior of $D_k$ near its special
points (tangencies, nodes, cusps). Moreover, one must take into
account the possible occurrence of creations or cancellations of admissible
pairs of nodes in the branch curve $D_k$, which affect the braid monodromy
morphism $\rho_k:\pi_1(\C-\{\mathrm{pts}\})\to B_d$ by insertion or deletion
of pairs of factors. The uniqueness statement in Theorem \ref{thm:cover}
then leads to the following result \cite{AK}:

\begin{theorem}[A.-Katzarkov]\label{thm:ak}
For given large enough $k$, the monodromy morphisms $(\rho_k,\theta_k)$
associated to the approximately holomorphic branched covering maps
$f_k:X\to\CP^2$ defined by triples of sections of $L^{\otimes k}$ are, up
to conjugation, Hurwitz operations, and insertions/deletions, invariants of
the symplectic manifold $(X,\omega)$. Moreover, these invariants are {\em
complete}, in the sense that the data $(\rho_k,\theta_k)$ are sufficient to
reconstruct the manifold $(X,\omega)$ up to symplectomorphism.
\end{theorem}

It is interesting to mention that the symplectic Lefschetz pencils
constructed by Donaldson (Theorem \ref{thm:donaldson}) can be recovered very easily from
the branched covering maps $f_k$, simply by considering the $\CP^1$-valued
maps $\pi\circ f_k$. In other words, the fibers $\Sigma_{k,\alpha}$ of the
pencil are the preimages by $f_k$ of the fibers of $\pi$, and the singular
fibers of the pencil correspond to the fibers of $\pi$ through the tangency
points of $D_k$.

In fact, the monodromy morphisms $\psi_k$ of the Lefschetz pencils $\pi\circ
f_k$ can be recovered very explicitly from $\theta_k$ and
$\rho_k$. By restriction to the line $\bar\ell=\ell\cup \{\infty\}$, the
$S_N$-valued morphism $\theta_k$ describes the topology of a fiber
$\Sigma_k$ of the pencil as an $N$-fold covering of $\CP^1$ with $d$ branch
points; the set of base points $Z_k$ is the preimage of the point at
infinity in $\bar\ell$. This makes it possible to define a {\it lifting
homomorphism} from a subgroup $B_d^0(\theta_k)\subset B_d$ ({\it liftable
braids}) to the mapping class group $\mathrm{Map}(\Sigma_k,Z_k)=
\mathrm{Map}_{g,N}$. The various monodromies are then related by the
following formula \cite{AK}:
\begin{equation}\label{eq:lifting}\psi_k=(\theta_k)_*\circ \rho_k.\end{equation}

The lifting homomorphism $(\theta_k)_*$ maps liftable half-twists to Dehn
twists, so that the tangencies between the branch curve $D_k$ and the fibers
of $\pi$ determine explicitly the vanishing cycles of the Lefschetz pencil
$\pi\circ f_k$. On the other hand, the monodromy around a node or cusp of
$D_k$ lies in the kernel of $(\theta_k)_*$.

The lifting homomorphism $\theta_*$ can be defined more precisely as
follows: the space $\tilde{\mathcal{X}}_d$ of configurations of $d$
distinct points in $\R^2$ together with branching data (a transposition
in $S_N$ attached to each point, or more accurately an $S_N$-valued group
homomorphism) is a finite covering of the space
$\mathcal{X}_d$ of configurations of $d$ distinct points. The morphism
$\theta$ determines a lift $\tilde{*}$ of the base point in $\mathcal{X}_d$,
and the liftable braid subgroup of $B_d=\pi_1(\mathcal{X}_d,*)$ is
the stabilizer of $\theta$ for the action of $B_d$ by deck transformations
of the covering $\tilde{\mathcal{X}}_d\to \mathcal{X}_d$, i.e.\
$B_d^0(\theta)=\pi_1(\tilde{\mathcal{X}}_d,\tilde{*})$.
Moreover, $\tilde{\mathcal{X}}_d$ is naturally equipped with a 
universal fibration $\mathcal{Y}_d\to \tilde{\mathcal{X}}_d$ by genus $g$
Riemann surfaces with $N$ marked points: the lifting homomorphism
$\theta_*:B_d^0(\theta)\to \mathrm{Map}_{g,N}$
is by definition the monodromy of this fibration.

The relation (\ref{eq:lifting}) is very useful for explicit calculations of
the monodromy of Lefschetz pencils, which is accessible to direct methods
only in a few very specific cases. By comparison, the various available
techniques for braid monodromy calculations \cite{Mo3,Te1,ADKY} are much
more powerful, and make it possible to carry out calculations in a much
larger number of cases (see below). In particular, in view of
Donaldson's result we are mostly interested in the monodromy of high degree
Lefschetz pencils, where the fiber genus and the number of singular fibers
are very high, making them inaccessible to direct calculation even for the
simplest complex algebraic surfaces.

When considering higher-dimensional manifolds, one may also associate monodromy
invariants to a fibration $f_k:X\to\CP^2$; these consist of the braid
monodromy of the critical curve $D_k\subset\CP^2$, and the monodromy of
the fibration (with values in the mapping class group of the fiber).
Furthermore, considering successive hyperplane sections and projections
to $\CP^2$, one obtains a complete description of a symplectic manifold
$(X^{2n},\omega)$ in terms of $n-1$ braid group-valued monodromy morphisms
(describing the critical curves of various $\CP^2$-valued projections)
and a single symmetric group-valued homomorphism (see \cite{Au3} for
details).
\medskip

In principle, the above results reduce the classification of compact
symplectic manifolds to purely combinatorial questions concerning braid
groups and symmetric groups, and symplectic topology seems to largely reduce
to the study of certain singular plane curves, or equivalently
factorizations in braid groups.

The explicit calculation of these monodromy invariants is hard in the
general case, but is made possible for a large number of complex surfaces by
the use of ``degeneration'' techniques and of approximately holomorphic
perturbations. These techniques make it possible to compute the braid
monodromies of a variety of algebraic branch curves. In most cases, the
calculation is only possible for a fixed projection to $\CP^2$ of a given
algebraic surface (i.e., fixing $[\omega]$ and considering only $k=1$);
see e.g.\ \cite{Mo2,Ro,AGTV,...} for various such examples. In a smaller
set of examples, the technique applies to projections of arbitrarily large
degrees and hence makes it possible to compute explicitly the invariants
defined by Theorem \ref{thm:ak}. To this date, the list consists of
$\CP^2$ \cite{MVeronese,TVeronese}, $\CP^1\times\CP^1$
\cite{Mo3}, certain Del Pezzo and K3 surfaces \cite{Ro,CMT},
the Hirzebruch surface $\mathbb{F}_1$ \cite{MRT,ADKY}, and all double covers
of $\CP^1\times\CP^1$ branched along connected smooth algebraic curves,
which includes an infinite family of surfaces of general type \cite{ADKY}.

The degeneration technique, developed by Moishezon and Teicher
\cite{Mo3,Te1,...}, starts with a projective embedding of the complex surface
$X$, and deforms the image of this embedding to a singular configuration
$X_0$ consisting of a union of planes intersecting along lines. The
discriminant curve of a projection of $X_0$ to $\CP^2$ is therefore a union
of lines; the manner in which the smoothing of $X_0$ affects this curve
can be studied explicitly, by considering a certain number of standard local
models near the various points of $X_0$ where three or more planes
intersect (see \cite{Te1} and \cite{ADKY} for detailed reviews of the
technique). This method makes it possible to handle many examples in low
degree, but in the case $k\gg 0$ that we are interested in (very positive
linear systems over a fixed manifold), the calculations can only
be carried out explicitly for very simple surfaces.

In order to proceed beyond this point, it becomes more efficient to move
outside of the algebraic framework and to consider generic approximately
holomorphic perturbations of non-generic algebraic maps; the greater
flexibility of this setup makes it possible to choose more easily computable
local models. For example, the direct calculation of the monodromy
invariants becomes possible for all linear systems of the type $\pi^*O(p,q)$
on double couvers of $\CP^1\times\CP^1$ branched along connected smooth
algebraic curves of arbitrary degree \cite{ADKY}. It also becomes possible
to obtain a general ``degree doubling'' formula, describing explicitly the
monodromy invariants associated to the linear system $L^{\otimes 2k}$ in
terms of those associated to the linear system $L^{\otimes k}$ (when $k\gg
0$), both for branched covering maps to $\CP^2$ and for 4-dimensional
Lefschetz pencils \cite{AK2}.

However, in spite of these successes, a serious obstacle restricts the
practical applications of monodromy invariants: in general, they cannot
be used efficiently to distinguish homeomorphic symplectic manifolds,
because no algorithm exists to decide whether two words in a braid group
or mapping class group are equivalent to each other via Hurwitz operations.
Even if an algorithm could be found, another difficulty is due to the large
amount of combinatorial data to be handled: on a typical interesting example,
the braid monodromy data can already consist of $\sim 10^4$ factors in
a braid group on $\sim 100$ strings for very small values of the
parameter $k$, and the amount of data grows polynomially with $k$.

Hence, even when monodromy invariants can be computed, they cannot be {\it
compared}. This theoretical limitation makes it necessary to search for
other ways to exploit monodromy data, e.g.\ by considering invariants that
contain less information than braid monodromy but are easier to use in
practice.

\subsection{Fundamental groups of branch curve complements}

Given a singular plane curve $D\subset\CP^2$, e.g.\ the branch curve of a
covering, it is natural to study the fundamental group $\pi_1(\CP^2-D)$.
The study of this group for various types of algebraic curves is a classical
subject going back to the work of Zariski, and has undergone a lot of
development in the 80's and 90's, in part thanks to the work of Moishezon
and Teicher \cite{Mo2,Mo3,Te1,...}. The relation to braid monodromy invariants
is a very direct one: the Zariski-van Kampen theorem provides an explicit
presentation of the group $\pi_1(\CP^2-D)$ in terms of the braid monodromy
morphism $\rho:\pi_1(\C-\{\mathrm{pts}\})\to B_d$. However, if one is
interested in the case of symplectic branch curves, it is important to
observe that the introduction or the cancellation of pairs of nodes affects
the fundamental group of the complement, so that it cannot be used directly
to define an invariant associated to a symplectic branched covering.
In the symplectic world, the fundamental group of the branch curve
complement must be replaced by a suitable quotient, the {\it stabilized
fundamental group} \cite{ADKY}.

Using the same notations as above, the inclusion \hbox{$i:\ell-(\ell\cap
D_k)\to \mathbb{CP}^2-D_k$} of the reference fiber of the linear projection
$\pi$ induces a surjective morphism on fundamental groups; the images of the
standard generators of the free group \hbox{$\pi_1(\ell-(\ell\cap D_k))$} and their
conjugates form a subset $\Gamma_k\subset\pi_1(\mathbb{CP}^2-D_k)$ whose
elements are called {\it geometric generators}. Recall that the images of
the geometric generators by the monodromy morphism $\theta_k$ are
transpositions in $S_N$. The creation of a pair of nodes in the curve $D_k$
amounts to quotienting $\pi_1(\CP^2-D_k)$ by a relation of the form
$[\gamma_1,\gamma_2]\sim 1$, where $\gamma_1,\gamma_2\in\Gamma_k$; however,
this creation of nodes can be carried out by deforming the branched
covering map $f_k$ only if the two transpositions $\theta_k(\gamma_1)$ and
$\theta_k(\gamma_2)$ have disjoint supports. Let $K_k$ be the normal
subgroup of $\pi_1(\CP^2-D_k)$ generated by all such commutators
$[\gamma_1,\gamma_2]$. Then we have the following result \cite{ADKY}:

\begin{theorem}[A.-D.-K.-Y.]
For given $k\gg 0$, the {\em stabilized fundamental group}
$\bar{G}_k=\pi_1(\mathbb{CP}^2-D_k)/K_k$ is an invariant of the symplectic
manifold $(X^4,\omega)$.
\end{theorem}

This invariant can be calculated explicitly for the various examples where
monodromy invariants are computable ($\CP^2$, $\CP^1\times\CP^1$, some
Del Pezzo and K3 surfaces, Hirzebruch surface $\mathbb{F}_1$,
double covers of $\CP^1\times\CP^1$); namely, the extremely complicated
presentations given by the Zariski-van Kampen theorem in terms of braid
monodromy data can be simplified in order to obtain a manageable
description of the fundamental group of the branch curve complement.
These examples lead to various observations and conjectures.

A first remark to be made is that, for all known examples, when the
parameter $k$ is sufficiently large the stabilization operation becomes
trivial, i.e.\ geometric generators associated to disjoint transpositions
already commute in $\pi_1(\CP^2-D_k)$, so that $K_k=\{1\}$ and
$\bar{G}_k=\pi_1(\CP^2-D_k)$. For example, in the case of $X=\CP^2$ with
its standard K\"ahler form, we have $\bar{G}_k=\pi_1(\CP^2-D_k)$ for all
$k\ge 3$. Therefore, when $k\gg 0$ no information seems to be lost when
quotienting by $K_k$ (the situation for small values of $k$ is very
different).

The following general structure result can be proved for the groups
$\bar{G}_k$ (and hence for $\pi_1(\CP^2-D_k)$) \cite{ADKY}:

\begin{theorem}[A.-D.-K.-Y.]
Let $f:(X,\omega)\to \CP^2$ be a symplectic branched covering of degree $N$,
with braided branch curve $D$ of degree $d$, and let
$\bar{G}=\pi_1(\CP^2-D)/K$ be the stabilized fundamental group of the branch
curve complement. Then there exists a natural exact sequence
$$1\longrightarrow G^0\longrightarrow \bar{G} \longrightarrow
S_N\times \Z_d\longrightarrow \Z_2\longrightarrow 1.$$
Moreover, if $X$ is simply connected
then there exists a natural surjective homomorphism
$\phi:G^0\twoheadrightarrow (\Z^2/\Lambda)^{n-1}$, where
$$\Lambda=\{(c_1(K_X)\cdot \alpha,[f^{-1}(\bar\ell)]\cdot \alpha),
\ \alpha\in H_2(X,\Z)\}.$$
\end{theorem}

In this statement, the two components of the morphism $\bar{G}\to S_N\times
\Z_d$ are respectively the monodromy of the branched covering,
$\theta:\pi_1(\CP^2-D)\to S_N$, and the linking number (or abelianization,
when $D$ is irreducible)
morphism $$\delta:\pi_1(\CP^2-D)\to \Z_d\ (\simeq H_1(\CP^2-D,\Z)).$$
The subgroup $\Lambda$ of $\Z^2$ is entirely determined by the numerical
properties of the canonical class $c_1(K_X)$ and of the hyperplane class
(the homology class of the preimage of a line $\bar\ell\subset\CP^2$: in
the case of the covering maps of Theorem \ref{thm:cover} we have 
$[f^{-1}(\bar\ell)]=c_1(L^{\otimes k})=\frac{k}{2\pi}[\omega]$). The morphism
$\phi$ is defined by considering the $N$ lifts in $X$ of a closed loop
$\gamma$ belonging to $G^0$, or more precisely their homology classes (whose
sum is trivial) in the complement of a hyperplane section and of the
ramification curve in $X$.

Moreover, in the known examples we have a much stronger result on
the structure of the subgroups $G_k^0$ for the branch
curves of large degree covering maps (determined by sufficiently ample
linear systems) \cite{ADKY}.

Say that the simply connected complex surface $(X,\omega)$ belongs to the
class $(\mathcal{C})$ if it belongs to the list of computable examples:
$\CP^1\times\CP^1$, $\CP^2$, the Hirzebruch surface $\mathbb{F}_1$ (equipped
with any K\"ahler form), a Del Pezzo or K3 surface (equipped with a K\"ahler
form coming from a specific complete intersection realization), or a double
cover of $\CP^1\times \CP^1$ branched along a connected smooth algebraic
curve (equipped with a K\"ahler form in the class $\pi^*O(p,q)$ for $p,q\ge
1$). Then we have:

\begin{theorem}[A.-D.-K.-Y.]
If $(X,\omega)$ belongs to the class $(\mathcal{C})$, then
for all large enough $k$ the homomorphism
$\phi_k$ induces an isomorphism on the abelianized groups, i.e.\ $\mathrm{Ab}
\,G_k^0\simeq$ $(\Z^2/\Lambda_k)^{N_k-1}$, while $\mathrm{Ker}\,\phi_k=[G_k^0,
G_k^0]$ is a quotient of $\Z_2\times\Z_2$.
\end{theorem}

It is natural to make the following conjecture:

\begin{conj}\label{conj:adky}
If $X$ is a simply connected symplectic 4-manifold, then for all large
enough $k$ the homomorphism
$\phi_k$ induces an isomorphism on the abelianized groups, i.e.\
$\mathrm{Ab} \,G_k^0\simeq$ $(\Z^2/\Lambda_k)^{N_k-1}$.
\end{conj}

\subsection{Symplectic isotopy and non-isotopy}

While it has been well-known for many years that compact symplectic
4-manifolds do not always admit K\"ahler structures, it has been discovered
more recently that symplectic curves (smooth or singular) in a given
manifold can also offer a wider range of possibilities than complex curves.
Proposition \ref{prop:cover} and Theorem \ref{thm:cover} establish a bridge between these two
phenomena: indeed, a covering of $\CP^2$ (or any other complex surface)
branched along a complex curve automatically inherits a complex structure.
Therefore, starting with a non-K\"ahler symplectic manifold, Theorem
\ref{thm:cover} always yields branch curves that are not isotopic to any complex curve in
$\CP^2$. The study of isotopy and non-isotopy phenomena for curves is
therefore of major interest for our understanding of the topology of
symplectic 4-manifolds.

The {\it symplectic isotopy problem} asks whether, in a given complex
surface, every symplectic submanifold representing a given homology class
is isotopic to a complex submanifold. The first positive result in this
direction was due to Gromov, who showed using his compactness result for
pseudo-holomorphic curves (Theorem \ref{thm:gromov}) that, in $\CP^2$, a smooth
symplectic curve of degree $1$ or $2$ is always isotopic to a complex curve.
Successive improvements of this technique have made it possible to extend
this result to curves of higher degree in $\CP^2$ or $\CP^1\times\CP^1$; the
currently best known result is due to Siebert and Tian, and makes it
possible to handle the case of smooth curves in $\CP^2$ up to degree 17
\cite{ST} (see also their lecture notes in this volume).
Isotopy results are also known for sufficiently simple singular
curves (Barraud, Shevchishin \cite{Sh}, \dots).

Contrarily to the above examples, the general answer to the symplectic
isotopy problem appears to be negative. The first counterexamples
among smooth connected symplectic curves were found by Fintushel and Stern
\cite{FS2}, who constructed by a {\it braiding} process infinite families of
mutually non-isotopic symplectic curves representing a same homology class
(a multiple of the fiber) in elliptic surfaces (a similar construction can
also be performed in higher genus). However, these two constructions are
preceded by a result of Moishezon \cite{Mo4}, who established in the
early 90's a result implying the existence in $\CP^2$ of infinite families
of pairwise non-isotopic singular symplectic curves of given degree with
given numbers of node and cusp singularities. A reformulation of Moishezon's
construction makes it possible to see that it also relies on braiding;
moreover, the braiding construction can be related to a surgery operation
along a Lagrangian torus in a symplectic 4-manifold, known as {\it Luttinger
surgery} \cite{ADK}. This reformulation makes it possible to vastly simplify
Moishezon's argument, which was based on lengthy and delicate calculations
of fundamental groups of curve complements, while relating it with various
constructions developed in 4-dimensional topology.

Given an embedded Lagrangian torus $T$ in a symplectic 4-manifold
$(X,\omega)$, a homotopically non-trivial embedded loop $\gamma\subset T$
and an integer $k$, Luttinger surgery is an operation that consists in
cutting out from $X$ a tubular neighborhood of $T$, foliated by parallel
Lagrangian tori, and gluing it back in such a way that the new meridian loop
differs from the old one by $k$ twists along the loop $\gamma$ (while
longitudes are not affected), yielding a new symplectic manifold
$(\tilde{X},\tilde{\omega})$. This relatively little-known construction,
which e.g.\ makes it possible to turn a product $T^2\times \Sigma$ into any
surface bundle over $T^2$, or to transform an untwisted fiber sum into a
twisted one, can be used to described in a unified manner numerous examples
of exotic symplectic 4-manifolds constructed in the past few years.

Meanwhile, the braiding construction of symplectic curves starts with a
(possibly singular) symplectic curve $\Sigma\subset(Y^4,\omega_Y)$ and two
symplectic cylinders embedded in $\Sigma$, joined by a
Lagrangian annulus contained in the complement of $\Sigma$,
and consists in performing $k$ half-twists between these two cylinders in
order to obtain a new symplectic curve $\tilde\Sigma$ in $Y$.

\begin{center}
\epsfig{file=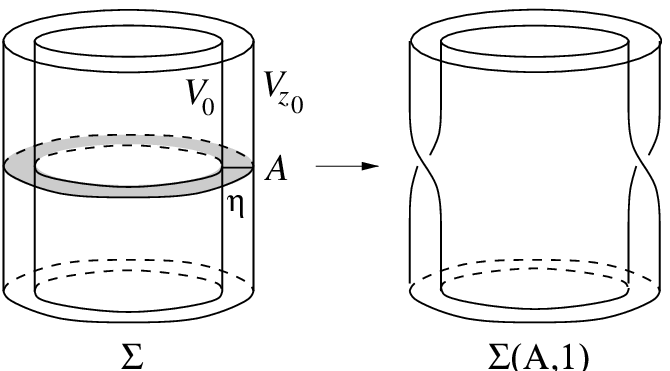,height=1.2in}
\end{center}

When $\Sigma$ is the branch curve of a symplectic branched covering $f:X\to
Y$, the following result holds \cite{ADK}:

\begin{proposition}\label{prop:adk}
The covering of $Y$ branched along the symplectic curve
$\tilde\Sigma$ obtained by braiding $\Sigma$ along a Lagrangian annulus
$A\subset Y-\Sigma$ is naturally symplectomorphic to the manifold
$\tilde{X}$ obtained from the branched cover $X$ by Luttinger surgery along
a Lagrangian torus $T\subset X$ formed by the union of two lifts of $A$.
\end{proposition}

Hence, once an infinite family of symplectic curves has been constructed by
braiding, it is sufficient to find invariants that distinguish the
corresponding branched covers in order to conclude that the curves are not
isotopic. In the Fintushel-Stern examples, the branched covers are
distinguished by their Seiberg-Witten invariants, whose behavior is well
understood in the case of elliptic fibrations and their surgeries.

In the case of Moishezon's examples, a braiding construction makes it
possible to construct, starting from complex curves
$\Sigma_{p,0}\subset\CP^2$ ($p\ge 2$) of degree $d_p=9p(p-1)$ with
$\kappa_p=27(p-1)(4p-5)$ cusps and $\nu_p=27(p-1)(p-2)(3p^2+3p-8)/2$ nodes,
symplectic curves $\Sigma_{p,k}\subset \CP^2$ for all $k\in\Z$, with the same
degree and numbers of singular points. By Proposition \ref{prop:adk}, these curves can
be viewed as the branch curves of symplectic coverings whose total spaces
$X_{p,k}$ differ by Luttinger surgeries along a Lagrangian torus $T\subset
X_{p,0}$. The effect of these surgeries on the canonical class and on the
symplectic form can be described explicitly, which makes it possible to
distinguish the manifolds $X_{p,k}$: the canonical class of
$(X_{p,k},\omega_{p,k})$ is given by
$p\,c_1(K_{p,k})=(6p-9)[\omega_{p,k}]+(2p-3)k\,PD([T])$. Moreover, 
$[T]\in H_2(X_{p,k},\Z)$ is not a torsion class, and if $p\not\equiv 0\mod
3$ or $k\equiv 0\mod 3$ then it is a primitive class \cite{ADK}. This implies
that infinitely many of the curves $\Sigma_{p,k}$ are pairwise non-isotopic.

It is to be observed that the argument used by Moishezon to distinguish the
curves $\Sigma_{p,k}$, which relies on a calculation of the
fundamental groups $\pi_1(\CP^2-\Sigma_{p,k})$ \cite{Mo4}, is related to
the one in \cite{ADK} by means of Conjecture \ref{conj:adky}, of which it can be
concluded a posteriori that it is satisfied by the given branched covers
$X_{p,k}\to \CP^2$: in particular, the fact that $\pi_1(\CP^2-\Sigma_{p,k})$
is infinite for $k=0$ and finite for $k\neq 0$ is consistent with the
observation that the canonical class of $X_{p,k}$ is proportional to its
symplectic class iff $k=0$.

\section{Symplectic surfaces from symmetric products}

This section will describe an approach, developed in \cite{DS},
\cite{Sm3} and \cite{Usher}, to various theorems due to Taubes,
but in the context of Lefschetz pencils rather than Seiberg-Witten gauge
theory and the equivalence ``SW=Gr''. The main result on which we will focus is Theorem \ref{thm:taubes}
$(i)$ about the existence of embedded symplectic submanifolds representing
the canonical class (although we will establish the result only in the case where $b_+>1+b_1$). Note that a Lefschetz pencil is a family of symplectic
surfaces of arbitrarily large volume and complexity (depending on $k$),
so somewhat different techniques are needed to find symplectic surfaces in
a {\it prescribed} homology class. 

\subsection{Symmetric products}

Let $(X,\omega_X)$ be a symplectic 4-manifold with
integral $[\omega_X]$, and fix a Lefschetz pencil, giving rise to a Lefschetz
fibration $f:\hat{X}\to S^2$ on the blow-up. We will always equip $\hat{X}$
with the symplectic form $\omega_C=p^*\omega_X+C\,f^*\omega_{S^2}$ for some
large $C>0$.

\begin{definition}
A standard surface in $\hat{X}$ is a smooth embedded surface $\Sigma$ such
that $f_{|\Sigma}$ is a simple branched covering over $S^2$, and at each
branch point $df$ gives an oriented isomorphism between the normal bundle
to $\Sigma$ and $TS^2$.
\end{definition}

The condition on $df$ ensures that we have positive branching, that is, near
each branch point there are complex-valued coordinates such that $\Sigma$
looks like $\{(z_1,z_2)\in\C^2,\ z_1^2=z_2\}$, with $f$ being the projection
to $z_2$. 

\begin{lemma}
A standard surface in $\hat{X}$, disjoint from the exceptional sections
$E_i$, defines a symplectic surface in $X$.
\end{lemma}

This follows from a local calculation showing that $(\hat{X}-\bigcup E_i,
\omega_C)$ is symplectomorphic to $(X-\{b_i\},(1+kC)\omega_X)$, where $k$
is the degree of the pencil. By taking $C$ sufficiently large it is clear
that the given standard surface is symplectic in $(\hat{X},\omega_C)$.

Rather than encoding a standard surface algebraically in the monodromy representation, we will take a geometric standpoint and obtain such surfaces from 
sections of a ``relative symmetric
product''.  This construction will in fact yield singular unions of such
surfaces, so for later use we therefore provide a smoothing lemma. As a piece
of notation, a {\it positive symplectic divisor} is a union of symplectic
surfaces $D=\sum a_i\Sigma_i$ with $a_i>0$, with the $\Sigma_i$ pairwise 
having only isolated positive transverse intersections, and no triple
intersections. (In fact, a careful inspection of the later argument will
show that, in the positive divisors we construct, if $\Sigma_i$ and
$\Sigma_j$ intersect, then either one is an exceptional sphere, or one
has multiplicity $a_i>1$.) In any case, we will be able to use the following
criterion:

\begin{proposition}
If a symplectic divisor satisfies $D\cdot\Sigma_j\ge 0$ for all $j$, then
it can be symplectically smoothed.
\end{proposition}

\proof[Sketch of proof]
Suppose first that all the $\Sigma_i$ are complex curves in a K\"ahler
surface. The hypothesis says that the line bundle $O(D)$ has non-negative
degree on each component. Then we can find a smooth section $\gamma$ of
$O(D)$ which is holomorphic near each intersection point and near its
(transverse) zero set. If the section $s$ defines the divisor $D$, then
the zero set of $s-\epsilon\gamma$ is a suitable symplectic smoothing for
sufficiently small $\epsilon$.

In general, the $\Sigma_i$ may not be symplectomorphic to any complex model.
Near an intersection point of $\Sigma$ ad $\Sigma'$, write $\Sigma'$ as a
graph of a linear map $A$ over the symplectic orthogonal of $\Sigma$. The
symplecticity of $\Sigma'$ implies that $\det(A)>-1$; for a complex model
(in which $A$ is an element of $\C\subset \mathrm{End}(\R^2)$), we need 
$\det(A)>0$ (or $A=0$). This can be achieved by a small compactly supported
perturbation of $\Sigma'$.
\endproof

We now turn to the main construction. The $d$-fold symmetric prouct of a
Riemann surface is canonically a smooth complex manifold of dimension $d$.
A holomorphic atlas is provided by noticing that the elementary symmetric
functions define a biholomorphism
$\mathrm{Sym}^d(\C)=\C^d/S_d\stackrel{\sim}{\longrightarrow}\C^d$.
For instance, taking a non-zero $(d+1)$-tuple of numbers to the roots of
the degree $d$ polynomial having these as its coefficients gives an
identification between $\CP^d$ and $\mathrm{Sym}^d(\CP^1)$.

Note that the symmetric product of a smooth two-manifold has no natural
smooth atlas, although it is well-defined up to non-canonical
diffeomorphisms. We will always work with almost-complex structures on our
Lefschetz fibrations which make the projection map pseudo-holomorphic (call
these almost-complex structures {\it fibered}). In
this case the fibers become Riemann surfaces.

\begin{theorem}
Given a Lefschetz fibration $f:\hat{X}\to\CP^1$ with a fibered
almost-complex structure, for each $d\ge 0$ there is a smooth compact
symplectic manifold $X_d(f)\to \CP^1$ with fiber at $t\not\in
\mathrm{crit}(f)$ the symmetric product $\mathrm{Sym}^d(f^{-1}(t))$.
\end{theorem}

Away from the singular fibers, this follows by taking charts on $X$ of
the form $\theta: D_1 \times D_2 \rightarrow X$, such that: for each $x \in
D_1$, $\theta_x: \{x\} \times D_2 \rightarrow X$ gives a holomorphic disc
in a fiber of $f$; and $f \circ \theta = i \circ \pi_1$, for $i$ an
inclusion $D_1 \hookrightarrow \CP^1$. In other words, the fibers of $f$
are Riemann surfaces whose complex structures vary smoothly in the obvious
way. Now we restrict further to work with almost complex structures which
are integrable near the singular fibers.  In this case, we can fill in the
fibration by appealing to algebraic geometry.  Namely, replace the symmetric product (a moduli space of
structure sheaves) by a \emph{Hilbert scheme} (a moduli space of ideal sheaves).
We will use a concrete description of the Hilbert scheme of $\C^2$ 
together with our local model for the singularities of $f$.

\begin{definition}
$Hilb_d(\C^2)$ comprises the triples $\{B_1,B_2,g\}/\sim$ where $B_1,B_2$
are $d \times d$ commuting matrices, and $g$ is a vector which is
not contained in any subspace $S \subset \C^d$ stabilized by both
$B_i$. $GL_d(\C)$ acts by conjugating the $B_i$'s and on the left of
$g$.  
\end{definition}

To unravel this, think of an ideal $J \subset \C[z_1,z_2]$ of codimension
$d$. Identify $\C[z_1,z_2]/J = \C^d$, and define $B_1,B_2$ to be the
action of $z_1,z_2$, and $g$ the image of $1$. Conversely, given a 
triple, define a map $\C[z_1,z_2] \rightarrow \C^d$ by $f(z_1,z_2)
\mapsto f(B_1,B_2)g$. The stability condition on $g$ implies that this is
surjective, and the kernel gives an ideal of codimension $d$.  Nakajima's
notes \cite{Nakajima} show that this correspondence is indeed a holomorphic
parametrisation of the Hilbert scheme of length $d$ subschemes of $\C^2$ by
triples of matrices.

Now consider the relative version of this, for the map $\C^2 \rightarrow
\C$, $(z_1,z_2) \mapsto z_1z_2$. The fiber over $t$ is by definition
those ideals which contain $z_1z_2-t$, giving 
\[
 B_1B_2 = t I, \quad B_2B_1 = t I, \quad \text{plus stability}.
\]
If the $B_i$ are simultaneously diagonalizable, we get a
$d$-tuple of pairs of eigenvalues, i.e.\ a point of $Sym^d(f^{-1}(t))$.
More generally, the $B_i$ can be represented as a pair of upper triangular
matrices, and again there is a map to the symmetric product given by taking
eigenvalues. For $t\neq 0$ this map is an isomorphism because of the stability
condition; for $t=0$ this is no longer the case, and suitable off-diagonal
matrix entries become coordinates on the fibers.  In any case, we now have explicit equations defining $X_d(f)$; one can use these to give a concrete description of its topology, and show smoothness.

\begin{proposition}
The total space of this relative Hilbert scheme is smooth.
\end{proposition}

\proof[Sketch of proof] One can show (for instance from the Abel-Jacobi discussion below)
that the $t = 0$ fiber  has normal crossings. Given this, at
each point of the singular locus there are two line bundles $\nu_i$
- the normal bundles to the two branches whose intersection defines the singular locus - and the given deformation of
the singular fiber gives a section of $\nu_1^* \otimes \nu_2^*$.
The total space is smooth if this section has no zeros. In our case,
the $\nu_i$ are canonically identified with the tangent spaces at the
node of the singular fiber of $f$, so the bundle $\nu_1^* \otimes \nu_2^*$ is trivial. Hence
we only need check smoothness at one point. This one can do by writing
down explicit local sections of the map.
\endproof

At least when we assume that all the singular fibers of $f$ are irreducible, we get another point of view from the Abel-Jacobi map. Let $\Sigma$
be a Riemann surface.

\begin{definition}
The Abel-Jacobi map $Sym^d\Sigma \rightarrow Pic^d(\Sigma)$ takes
a divisor $D$ to the line bundle $\mathcal{O}(D)$.
\end{definition}

The key point is to identify the fibers of this map with linear
systems $\mathbb{P}(H^0(L))$, hence projective spaces whose dimension may
vary as $L$ moves in $Pic$. 
\medskip

{\bf Example.} Let $d = 2g-2$ (by the adjunction formula, this is the relevant case if a standard surface in $\hat{X}$ is to represent $K_{\hat{X}}$). From the Riemann-Roch theorem $h^0(L) -
h^1(L) = d-g+1 = g-1$, whilst by Serre duality, $h^1(L) > 0$ only if
$L = K$ (since only the trivial degree $0$ bundle has a section).
Thus $Sym^{2g-2}(\Sigma) \rightarrow Pic^{2g-2}(\Sigma)$ is a
(locally trivial) $\CP^{g-2}$-bundle away from a single 
exceptional fiber, which is a $\CP^{g-1}$. 
For instance, $Sym^2(\Sigma_2) \rightarrow T^4$ blows down a single
exceptional sphere.
\medskip

A line bundle on a nodal Riemann surface $\Sigma_0$ is given by a line bundle
on its normalization $\tilde\Sigma_0$ together with a $\C^*$ gluing parameter 
$\lambda$ to identify the fibers over the preimages of the node. This 
$\C^*$-bundle over $T^{2g-2}$ naturally compactifies to a $\CP^1$-bundle;
if we glue the $0$ and $\infty$ sections of this over a translation
in the base, we compactify the family of Picard varieties for a
family of curves with an irreducible nodal fiber. (Think of an
elliptic Lefschetz fibration with section, which does precisely this
in the $g=1$ case!)  Thus, under our assumptions, the relative Picard
fibration $P_d(f)$ is also a smooth compact symplectic manifold.

Given $L \in Pic^{2g-2}(\Sigma_0)$ in the smooth locus, the line 
bundle on the normalization has a space of sections of dimension
$g$. Hence the $\lambda$-hyper\-plane of sections which transform
correctly over the node gives a $\CP^{g-2}$. The situation is 
similar along the normal crossing divisor in $Pic$ (just take
$\lambda = 0$ or $\infty$), but changes if $L$ gives rise to the
canonical bundle of the normalization, where the dimension jumps. Summing up:

\begin{proposition} \label{ajmap}
The symplectic manifold $X_{2g-2}(f)$ is the total space of a family of projective spaces over $P_{2g-2}(f)$, the family being a locally trivial $\CP^{g-1}$ bundle over the section defined by the canonical line bundles of the
fibers of $f$ and a locally trivial $\CP^{g-2}$ bundle away from this section.
\end{proposition}

There is an analogous description for other values of $d$, and one sees that the singular locus of the singular
fiber of $X_d(f)$ is a copy of $Sym^{d-1}(\tilde{\Sigma}_0)$.  It may be
helpful to point out that tuples set-theoretically including the node
are {\it not} the same thing as singular points of $Hilb(\Sigma_0)$. The
relevance is that a smooth family of surfaces in $X$ may well include members
passing through critical points of $f$, but smooth sections of $X_d(f)$ can
never pass through the singular loci of fibers.  (Roughly speaking, the singular locus of $Hilb_d(\Sigma_0)$ comprises tuples which contain the node with odd multiplicity.) 

To see that $X_d(f)$ admits symplectic structures, one uses an analogue of the argument that applies to
Lefschetz fibrations -- we have a family of K{\"a}hler manifolds
with locally holomorphic singularities, so we can patch together local forms
to obtain something vertically non-degenerate and then add the pull-back
of a form from the base. 

\subsection{Taubes' theorem}

The construction of symplectic surfaces representing $[K_X]$ proceeds in two
stages.
Let $r=2g-2$.
A section $s$ of $X_r(f)$ defines a cycle $C_s \subset \hat{X}$. There is a homotopy
class $h$ of sections such that $[C_h] = PD(-c_1(\hat{X}))$.
First, we show that the Gromov invariant counting sections of $X_r(f)$ in the class $h$ is non-zero, and then from the pseudoholomorphic sections of $X_r(f)$ thereby provided we construct standard surfaces (which have a suitable component which descends from $\hat{X}$ to $X$).  For the first stage the Abel-Jacobi map is the key, whilst for the second we work with almost complex structures compatible with the ``diagonal'' strata of the symmetric products (tuples of points not all of which are pairwise distinct).

That the Gromov invariant
is well-defined follows easily from a description of the possible bubbles arising from cusp curves in $X_r(f)$. Again, we will always use "fibered" almost complex
structures $J$ on $X_r(f)$, i.e. ones making the projection to $\CP^1$
holomorphic.  

\begin{lemma}
$\pi_2(Sym^r(\Sigma)) = \Z$, generated by a line $l$ in a fiber
of the Abel-Jacobi map. All bubbles in cusp limits of sections
of $X_r(f) \rightarrow \CP^1$ are homologous to a multiple
of $l$.
\end{lemma}

\proof[Sketch of proof] The first statement follows from the Lefschetz hyperplane theorem and the fact that
$Sym^d(\Sigma)$ is a projective bundle for $d > r$; it is also an easy consequence of the Dold-Thom theorem in algebraic topology.
The second statement is then obvious away from singular fibers, but these require special treatment, cf. \cite{DS}.
\endproof

Since $\langle c_1(Sym^r(\Sigma)), h \rangle > 0$ this shows that spaces of holomorphic sections can be compactified by high codimension pieces, after which it is easy to define the Gromov invariant.  

\begin{proposition}
For a suitable $J$, the
moduli space of $J$-holomorphic sections in the class $h$ is
a projective space of dimension $(b_+ - b_1 - 1)/2 - 1$. 
\end{proposition}

\proof The index of the $\bar\partial$-operator on
$f_*K \rightarrow \CP^1$ is $(b_+ - b_1 - 1)/2$.  Choose $J$ on $X_r(f)$ so
that the projection $\tau$ to $P_r(f)$ is holomorphic,
and the latter has the canonical section $s_K$ as an isolated 
holomorphic section (the index of the normal bundle to this
is $-(b_+ - b_1 -1)/2$). Now all holomorphic sections of $X_r(f)$ lie over $s_K$.  If the degree of the original Lefschetz pencil was sufficiently large, the rank of the bundle $f_*K$ is much larger than its
first Chern class. Then, for a generic connection, the bundle
becomes $\mathcal{O} \oplus \mathcal{O} \oplus \dots
\oplus \mathcal{O}(-1) \oplus \dots \oplus \mathcal{O}(-1)$
(by Grothendieck's theorem). Hence, all nonzero sections
are nowhere zero, so yield homotopic sections of the projectivisation
$\mathbb{P}(f_*K) = \tau^{-1}(s_K)$, each defining a cycle in $\hat{X}$ in the fixed 
homology class $-c_1(\hat{X})$. 
\endproof

We now have a $J$ such that $\mathcal{M}_J(h)$ is compact and smooth --
but of the wrong dimension (index arguments show that
for a section $s$ of $X_r(f)$ giving a cycle $C_s$ in $X$,
the space of holomorphic sections has virtual dimension 
$C_s \cdot C_s - C_s \cdot K_X$). The actual invariant is
given as the Euler class of an {\em obstruction bundle}, whose fiber at a
curve $C$ is $H^1(\nu_C)$ (from the above, all our curves are embedded
so this bundle is well-defined).  Let $Q$ be the {\it quotient bundle} over
$\CP^n$ (defined as the cokernel of the inclusion of the tautological
bundle into $\underline{\C}^{n+1}$). It has Euler class $(-1)^{n+1}$.

\begin{lemma}
For $J$ as above, the obstruction bundle $Obs \rightarrow \mathcal{M}_J(h)$ is
\,$Q\rightarrow \CP^N$ with $N = (b_+ - b_1 - 1)/2-1$.
\end{lemma}

\proof[Sketch of proof] It is not hard to show that there is a global
holomorphic model $M(f_*K)$ for the map $\tau: X_r(f)\rightarrow P_r(f)$
in a neighborhood of the fibers of excess dimension, i.e.
$\tau^{-1}(s_K) = \mathbb{P}(f_*K)$.  In other words, the obstruction
computation reduces to algebraic geometry (hence our notation for the
bundle $f_* K \rightarrow \CP^1$ with fiber $H^0(K_{\Sigma})$). We have a sequence of holomorphic vector bundles
$$
0 \rightarrow T(\mathbb{P}(f_*K)) \rightarrow
TM(f_*K) \rightarrow T(W) \rightarrow \mathrm{coker}(D\tau) 
\rightarrow 0 
$$
where $W$ is fiberwise dual to $f_*K$, but globally twisted
by $\mathcal{O}(-2)$ (this is because $K_{\hat{X}}|f^{-1}(t)$ is not canonically $K_{f^{-1}(t)}$, introducing a twist into all the
identifications). If $\tau$ was actually a submersion, clearly the obstruction bundle would be trivial (as $\nu_C$ would be
pulled back from $\tau(C)=s_K \subset P_r(f)$); the deviation from this is measured by 
$H^1(\mathrm{coker}(D\tau))$. We compute this by splitting the above sequence into two short exact sequences and
take the long exact sequences in cohomology.
\endproof

Summing up, when $b_+ > 1+b_1(X)$, we have shown that the Gromov invariant for the class $h$ is $\pm 1$.  (One can improve the constraint to $b_+>2$, but curiously it seems tricky to exactly reproduce Taubes' sharp bound $b_+ \geq 2$.)

The $J$-holomorphic curves in $X_r(f)$ give cycles in $\hat{X}$,
but there is no reason these should be symplectic. We are
trying to build standard surfaces, i.e.\ we need our 
cycles to have positive simple tangencies to fibers. This is
exactly saying we want sections of $X_r(f)$ \emph{transverse} to the
diagonal locus $\Delta$ and intersecting it \emph{positively}.
To achieve this, we will need to construct almost complex structures which behave well with respect to $\Delta$, which has a natural stratification by the combinatorial type of the tuple. Problems arise because we cannot a priori suppose that our sections do not have image entirely contained inside the diagonals.

As a first delicacy,
the diagonal strata are not smoothly embedded, but are all the
images of smooth maps which are generically bijective.  This leads to ``smooth models'' of the diagonal strata, and for suitable $J$ any holomorphic section will lift to a unique such model.  (In particular, if a given section lies inside $\Delta$ because it contains multiple components, then we ``throw away'' this multiplicity; restoring it later will mean that we construct positive symplectic divisors, in the first instance, and not embedded submanifolds.)  As a second delicacy, we haven't controlled the index of sections once we
lift them to smooth models of the strata.  Finally, although we are restricting to the class of almost complex structures which are compatible with $\Delta$, nonetheless we hope to find holomorphic sections transverse to $\Delta$ itself.  The result we need, then, can loosely be summed up as:

\begin{proposition}
There are ``enough'' almost complex structures compatible with
the strata.
\end{proposition}

A careful proof is given in \cite{DS}, Sections 6--7. The 
key is that we can obtain deformations of such almost complex structures
from vector fields on the fibers of $f$. Since we can find a
vector field on $\C$ taking any prescribed values at a fixed
set of points, once we restrict to an open subset of the graph of any given
holomorphic section which lies in the top stratum of its associated model
there are enough deformations to generate the
whole tangent space to the space of $J$ which are both fibered and compatible with $\Delta$.


\medskip

{\bf Remark.} There are also strata coming from exceptional sections
inside $X$, each with fiber $Sym^{2g-3}$. These also vary
topologically trivially, and we repeat the above discussion
for them.  This lets us get standard surfaces in $\hat{X} \backslash \cup E_i \cong X \backslash \{b_j\}$.
\medskip

The upshot is that the nonvanishing of the Gromov invariant for the class $h \in H_2(X_r(f))$ means we can find
a holomorphic section of $X_r(f)$ which is transverse and
positive to all strata of the diagonal in which it is not
contained.  Taking into account the strata coming from exceptional sections as well, and with care, one obtains  a positive symplectic divisor $D=\sum a_i\Sigma_i$ in $X$ in the homology class dual to $K_X$.  We must finally check that we can satisfy the conditions of our earlier smoothing result.
By adjunction, 
\[
 D \cdot \Sigma_j + \Sigma_j^2 = 2g(\Sigma_j) - 2 = (1 + a_j)\Sigma_j^2 +
 \sum_{i \neq j} a_i \, \Sigma_i \cdot \Sigma_j
\]
For $D\cdot \Sigma_j$ to be negative, it must be that $\Sigma_j$ is a 
rational $-1$ curve; but clearly it suffices to prove the theorem under the assumption that $X$ is minimal. Provided $b_+>1+b_1(X)$, we therefore obtain an embedded symplectic
surface in $X$ representing the class $K_X$.  (Strictly we have proven this theorem only when $\omega_X$ is rational, but one can remove this assumption, cf. \cite{Sm3}.)  

This result gives ``gauge-theory free'' proofs of various standard facts on symplectic four-manifolds: for instance, 
if $b_+(X) > 1 + b_1(X)$ and $X$ is minimal,
then $c_1^2(X) \geq 0$.  In particular, manifolds such as
$K3 \# K3 \# K3$ admit no symplectic structure, illustrating the claim (made in the Introduction) that the symplectic condition cannot be reduced to cohomological or homotopical conditions, in contrast to the case of open manifolds covered by Gromov's h-principle.

To close, we point out that 
we can also see Taubes' remarkable duality $Gr(D) = \pm
Gr(K-D)$ very geometrically in this picture. Let
$\iota: Pic^{2g-2-r}(\Sigma) \rightarrow Pic^r(\Sigma)$ denote the map $O(D) \mapsto O(K-D)$ and write $\tau_d$ for the Abel-Jacobi map $Sym^d(\Sigma) \rightarrow Pic^d(\Sigma)$.  If 
$r$ is small relative to $g$, 
then generically we expect the map $\tau_r$ to be an
embedding, and its image to be exactly where the map $\iota \circ \tau_{2g-2-r}$ has fibers of excess dimension. 
If we fix a homology class on $X$ and take pencils of
higher and higher degree, we can make $r$ (the intersection number
with the fiber) as small as we like relative to $g$ (the genus of the
fiber) -- the
latter grows quadratically, not linearly.  Combining this with some deep results in the Brill-Noether theory of Riemann surfaces, one can show the following generalisation of Proposition \ref{ajmap}:

\begin{proposition}[\cite{Sm3}]
For $r$ small enough relative to $g$, there exists $J$ on
$X_{2g-2-r}(f)$ such that $X_r(f)$ is embedded in
$P_r(f) \cong P_{2g-2-r}(f)$, and such that $X_{2g-2-r}(f) 
\rightarrow P_{2g-2-r}(f)$ is holomorphic
and locally trivial over $X_r(f)$ and its
complement.
\end{proposition} 

Using this and repeating all the above, we can show there are well-defined
integer-valued invariants $\mathcal{I}_f$ counting sections of the $X_r(f)$'s which satisfy
the duality 
$$\mathcal{I}_f(D) = \pm \mathcal{I}_f(K_X-D).$$

In fact, a smooth fibered $J$ on $X$ defines
a canonical fibered $\mathbb{J}$ on $X_r(f)$ but this is only $C^0$ along
the diagonal strata. Suitably interpreted, $J$-curves in
$X$ and $\mathbb{J}$-curves in $X_r(f)$ are tautologically 
equivalent. Motivated by this, and a (very rough) sketch showing 
$\mathcal{I}_f = Gr \!\mod 2$ for symplectic manifolds with $K_X = \lambda[\omega]$ for non-zero $\lambda$ (more generally whenever there are no embedded square zero symplectic tori), \cite{Sm3} conjectured that the $\mathcal{I}$-invariants and Taubes' Gromov invariants co-incide.  Recent work of Michael Usher \cite{Usher} clarifies and completes this circle of ideas by showing that in full generality $Gr=\mathcal{I}_f$ (in $\Z$) for any pencil $f$ of sufficiently high degree.  In this sense, the above equation shows that Taubes' duality can be understood in terms of Serre duality on the fibers of a Lefschetz fibration.

\section{Fukaya categories and Lefschetz fibrations}

Many central questions in symplectic topology revolve around Lagrangian
submanifolds, their existence and their intersection properties.
From a formal point of view
these can be encoded into so-called Fukaya categories, which also play
an important role of their own in Kontsevich's homological mirror symmetry conjecture.
It turns out that Lefschetz pencils provide some of the most powerful
tools available in this context, if one considers vanishing cycles
and the ``matching paths'' between them \cite{Se3}; the most exciting
applications to this date are related to \nobreak{Seidel's} construction
of ``directed Fukaya categories'' of vanishing cycles
\cite{Se1,Se2}, and to verifications of the homological mirror conjecture \cite{Se4}.

Fukaya categories of symplectic manifolds are intrinsically very hard to
compute, because relatively little is known about embedded Lagrangian
submanifolds in symplectic manifolds of dimension $4$ or more, especially in
comparison to the much better understood theory of coherent sheaves over
complex varieties, which play the role of their mirror counterparts.
The input provided by Lefschetz fibrations is a reduction
from the symplectic geometry of the total space to that of the fiber, which
in the 4-dimensional case provides a crucial simplification.

\subsection{Matching paths and Lagrangian spheres}
Let $f:X\to S^2$ be a symplectic Lefschetz fibration
(e.g.\ obtained by blowing up the base points of a Lefschetz pencil or,
allowing the fibers to be open, by simply removing them), and let
$\gamma\subset S^2$ be an embedded arc joining a regular value $p_0$
to a critical value $p_1$, avoiding all the other critical values of $f$.
Using the horizontal distribution given by the symplectic orthogonal to
the fibers, we can transport the vanishing cycle at $p_1$ along the arc
$\gamma$ to obtain a Lagrangian disc $D_\gamma\subset X$ fibered above $\gamma$,
whose boundary is an embedded Lagrangian sphere $S_\gamma$ in the fiber
$\Sigma_0=f^{-1}(p_0)$. The Lagrangian disc $D_\gamma$ is called
the {\it Lefschetz thimble} over $\gamma$, and its boundary $S_\gamma$ is the
vanishing cycle already considered in Section 2.

If we now consider an
arc $\gamma$ joining two critical values $p_1,p_2$ of $f$ and passing
through $p_0$, then the above construction applied to each half of $\gamma$
yields two Lefschetz thimbles $D_1$ and $D_2$, whose boundaries are
Lagrangian spheres $S_1,S_2\subset \Sigma_0$. If $S_1$ and $S_2$ coincide
exactly, then $D_1\cup D_2$ is an embedded Lagrangian sphere in $X$,
fibering above the arc $\gamma$ (see the picture below); more generally,
if $S_1$ and $S_2$ are Hamiltonian isotopic to each other, then perturbing
slightly the symplectic structure we can reduce to the previous case and
obtain again a Lagrangian sphere in $X$. The arc
$\gamma$ is called a {\it matching path} in the Lefschetz fibration $f$
\cite{Se3}.

\begin{figure}[ht]
\centerline{
\setlength{\unitlength}{0.5mm}
\begin{picture}(50,55)(-10,-15)
\put(-10,0){\line(1,0){50}}
\put(-10,40){\line(1,0){50}}
\qbezier[30](30,40)(35,40)(35,35)
\qbezier[30](35,35)(35,30)(34,27.5)
\qbezier[30](34,27.5)(33,25)(33,20)
\qbezier[30](30,0)(35,0)(35,5)
\qbezier[30](35,5)(35,10)(34,12.5)
\qbezier[30](34,12.5)(33,15)(33,20)
\qbezier[30](30,40)(25,40)(25,36)
\qbezier[30](25,36)(25,32)(26,28)
\qbezier[20](26,28)(27,26)(27,25)
\qbezier[30](30,0)(25,0)(25,5)
\qbezier[30](25,5)(25,10)(26,12.5)
\qbezier[80](26,12.5)(27,15)(27,25)
\qbezier[30](32,32)(32,36)(30,36)
\qbezier[30](30,36)(28,36)(26,32)
\qbezier[20](32,32)(32,30)(30,30)
\qbezier[20](30,30)(28,30)(26,32)
\qbezier[30](29,5)(33,8)(29,11)
\qbezier[20](30,6)(28,8)(30,10)
\qbezier[30](0,40)(5,40)(5,35)
\qbezier[30](5,35)(5,30)(4,27.5)
\qbezier[30](4,27.5)(3,25)(3,20)
\qbezier[30](0,0)(5,0)(5,5)
\qbezier[30](5,5)(5,10)(4,12.5)
\qbezier[30](4,12.5)(3,15)(3,20)
\qbezier[30](0,40)(-5,40)(-5,36)
\qbezier[30](-5,36)(-5,32)(-4,28)
\qbezier[20](-4,28)(-3,26)(-3,25)
\qbezier[30](0,0)(-5,0)(-5,5)
\qbezier[30](-5,5)(-5,10)(-4,12.5)
\qbezier[80](-4,12.5)(-3,15)(-3,25)
\qbezier[30](2,32)(2,36)(0,36)
\qbezier[30](0,36)(-2,36)(-4,32)
\qbezier[20](2,32)(2,30)(0,30)
\qbezier[20](0,30)(-2,30)(-4,32)
\qbezier[30](-1,5)(3,8)(-1,11)
\qbezier[20](0,6)(-2,8)(0,10)
\put(-4,32){\circle*{1.5}}
\put(26,32){\circle*{1.5}}
\qbezier[80](-4,32)(-4,36)(11,36)
\qbezier[80](-4,32)(-4,28)(11,28)
\qbezier[80](26,32)(26,28)(11,28)
\qbezier[80](26,32)(26,36)(11,36)
\qbezier[30](11,28)(9,32)(11,36)
\qbezier[10](11,28)(13,32)(11,36)
\put(9.5,20){$S^n$}
\put(-12,-15){\line(1,0){44}}
\put(-12,-15){\line(1,1){8}}
\put(32,-15){\line(1,1){8}}
\put(0,-10){\circle*{1.5}}
\put(30,-10){\circle*{1.5}}
\put(0,-10){\line(1,0){30}}
\put(15,-7){$\gamma$}
\end{picture}
}
\end{figure}

Matching paths are an important source of Lagrangian spheres, and more
generally (extending suitably the notion of matching path to embedded
arcs passing through several critical values of $f$) of embedded
Lagrangian submanifolds. Conversely, a folklore theorem asserts that
any given embedded Lagrangian sphere 
in a compact symplectic manifold is isotopic to one that fibers
above a matching path in a Donaldson-type symplectic Lefschetz pencil of
sufficiently high degree. More precisely, the following result holds
(see \cite{AMP} for a detailed proof):

\begin{theorem}\label{thm:matching}
Let $L$ be a compact Lagrangian submanifold in a compact symplectic manifold
$(X,\omega)$ with integral $[\omega]$, and let $h:L\to [0,1]$ be any Morse
function. Then for large enough $k$ there exist Donaldson pencils
$f_k:X-\{\mbox{base points}\}\to S^2$, embedded arcs $\gamma_k:[0,1]\hookrightarrow S^2$,
and Morse functions $h_k:L\to [0,1]$ isotopic to $h$,
such that the restriction of $f_k$ to $L$ is equal to $\gamma_k\circ h_k$.
\end{theorem}

The intersection theory of Lagrangian spheres that fiber above matching
paths is much nicer than that of arbitrary Lagrangian spheres, because
if two Lagrangian spheres $S,S'\subset X$ fiber above matching paths
$\gamma,\gamma'$, then all intersections of $S$ with $S'$ lie in the
fibers above the intersection points of $\gamma$ with $\gamma'$. Hence,
the Floer homology of $S$ and $S'$ can be computed by studying intersection
theory for Lagrangian spheres in the fibers of $f$ rather than in $X$.
In particular, when $X$ is a four-manifold the vanishing cycles are just
closed loops in Riemann surfaces, and the computation of Floer homology
essentially reduces to a combinatorial count.

The enumeration of matching paths, if possible, would lead to a complete
understanding of isotopy classes of Lagrangian spheres in a given symplectic
manifold, with various topological consequences and applications to the
definition of new symplectic invariants. However, no finite-time algorithm
is currently available for this problem, although some
improvements on the ``naive'' search are possible (e.g.\ using maps to
$\CP^2$ and projective duality to identify certain types of pencil
automorphisms). Nonetheless, Theorem \ref{thm:matching} roughly says that
all Lagrangians are built out of Lefschetz thimbles. This implies that,
at the formal level of (derived) Fukaya categories, it is sometimes possible
to identify a ``generating'' collection of Lagrangian submanifolds (out
of which all others can be built by gluing operations, or, in the language
of categories, by passing to bounded complexes). A spectacular illustration
is provided by Seidel's recent verification of homological mirror symmetry
for quartic K3 surfaces \cite{Se4}.

\subsection{Fukaya categories of vanishing cycles}

The above considerations, together with ideas of Kontsevich about
mirror symmetry for Fano varieties, have led
\nobreak{Seidel} to the following construction of a
Fukaya-type $A_\infty$-category associated to a symplectic Lefschetz
pencil $f$ on a compact symplectic manifold $(X,\omega)$ \cite{Se1}.
Let $f$ be a symplectic Lefschetz pencil determined by two sections
$s_0,s_1$ of a sufficiently positive line bundle $L^{\otimes k}$ as in
Theorem \ref{thm:donaldson}.
Assume that $\Sigma_\infty=s_1^{-1}(0)$ is a smooth fiber of the pencil,
and consider the symplectic manifold with boundary $X^0$ obtained from $X$
by removing a suitable neighborhood of $\Sigma_\infty$. The map $f$ induces
a Lefschetz fibration $f^0:X^0\to D^2$ over a disc,
whose fibers are symplectic submanifolds with boundary obtained from the
fibers of $f$ by removing a neighborhood of their intersection points with
the symplectic hypersurface $\Sigma_\infty$ (the base points of the pencil).
Choose a reference point $p_0\in \partial D^2$, and consider the fiber
$\Sigma_0=(f^0)^{-1}(p_0)\subset X^0$.

Let $\gamma_1,\dots,\gamma_r$ be a collection of arcs in $D^2$ joining the
reference point $p_0$ to the various critical values of $f^0$, intersecting
each other only at $p_0$, and ordered in the clockwise direction around
$p_0$. As discussed above, each arc $\gamma_i$ gives rise to a Lefschetz
thimble $D_i\subset X^0$, whose boundary is a Lagrangian sphere $L_i\subset
\Sigma_0$. To avoid having to discuss the orientation of moduli spaces, we
give the following definition using $\Z_2$ (instead of $\Z$) as the
coefficient ring \cite{Se1}:

\begin{definition}[Seidel]\label{def:seidel}
The Fukaya category of vanishing cycles $\mathcal{F}_{vc}(f;\{\gamma_i\})$ is a
(directed) $A_\infty$-category with $r$ objects $L_1,\dots,L_r$ (corresponding to
the vanishing cycles, or more accurately to the thimbles); the morphisms between
the objects are given by
$$\mathrm{Hom}(L_i,L_j)=\begin{cases}
CF^*(L_i,L_j;\Z_2)=\Z_2^{|L_i\cap L_j|} & \mathrm{if}\ i<j\\
\Z_2\,\,id & \mathrm{if}\ i=j\\
0 & \mathrm{if}\ i>j;
\end{cases}$$
and the differential $\mu^1$, composition $\mu^2$ and higher order
compositions $\mu^n$ are given by Lagrangian Floer homology inside
$\Sigma_0$. More precisely,
$$\mu^n:\mathrm{Hom}(L_{i_0},L_{i_1})\otimes \dots\otimes
\mathrm{Hom}(L_{i_{n-1}},L_{i_n}) \to \mathrm{Hom}(L_{i_0},L_{i_n})[2-n]$$
is trivial when the inequality $i_0<i_1<\dots<i_n$ fails to hold (i.e.\ it
is always zero in this case, except for $\mu^2$ where composition with an
identity morphism is given by the obvious formula).
When $i_0<\dots<i_n$, $\mu^n$ is defined
by counting pseudo-holomorphic maps from the disc to $\Sigma_0$, mapping
$n+1$ cyclically ordered marked points on the boundary to the given
intersection points between vanishing cycles, and the portions of boundary
between them to $L_{i_0},\dots,L_{i_n}$.
\end{definition}

One of the most attractive features of this definition is that it only
involves Floer homology for Lagrangians inside the hypersurface $\Sigma_0$;
in particular, when $X$ is a symplectic 4-manifold, the definition becomes
purely combinatorial, since in the case of a Riemann surface
the pseudo-holomorphic discs appearing in the definition of Floer homology
and product structures are just immersed polygonal regions with convex
corners.

  From a technical point of view, a property that greatly facilitates the
definition of Floer homology for the vanishing cycles $L_i$ is {\it
exactness}. Namely, the symplectic structure on the manifold
$X^0$ is {\it exact}, i.e.\ it can
be expressed as $\omega=d\theta$ for some 1-form $\theta$ (up to a scaling
factor, $\theta$ is the 1-form describing the connection on $L^{\otimes k}$
in the trivialization of $L^{\otimes k}$ over $X-\Sigma_\infty$ induced by
the section $s_1/|s_1|$). With this understood, the submanifolds
$L_i$ are all {\it exact Lagrangian}, i.e.\ the restriction $\theta_{|L_i}$
is not only closed ($d\theta_{|L_i}=\omega_{|L_i}=0$) but also exact,
$\theta_{|L_i}=d\phi_i$. Exactness has two particularly nice consequences.
First, $\Sigma^0$ contains no closed pseudo-holomorphic curves (because
the cohomology class of $\omega=d\theta$ vanishes). Secondly, there are no
non-trivial pseudo-holomorphic discs in $\Sigma_0$ with boundary contained
in one of the Lagrangian submanifolds $L_i$.
Indeed, for any such disc $D$, we have $\mathrm{Area}(D)=\int_D\omega=
\int_{\partial D}\theta=\int_{\partial D} d\phi_i=0$. Therefore, bubbling
never occurs (neither in the interior nor on the boundary of the domain) in
the moduli spaces used to define the Floer homology groups
$HF(L_i,L_j)$. Moreover, the exactness of $L_i$ provides a priori estimates
on the area of all pseudo-holomorphic discs contributing to the definition
of the products $\mu^n$ ($n\ge 1$); this implies the
finiteness of the number of discs to be considered and solves elegantly the
convergence problems that normally make it necessary to define Floer
homology over Novikov rings. 

\medskip

{\bf Example.} We now illustrate the above definition by considering
the example of a pencil of degree $2$ curves in $\CP^2$ (see also \cite{Se2}).
Consider the two sections
$s_0=x_0(x_1-x_2)$ and $s_1=x_1(x_2-x_0)$ of the line bundle $O(2)$ over
$\CP^2$: their zero sets are singular conics, in fact the unions of two
lines each containing two of the four intersection points $(1\!:\!0\!:\!0)$,
$(0\!:\!1\!:\!0)$, $(0\!:\!0\!:\!1)$, $(1\!:\!1\!:\!1)$. Moreover, the
zero set of the linear combination $s_0+s_1=x_2(x_1-x_0)$ is also singular;
on the other hand, it is fairly easy to check that all other linear
combinations $s_0+\alpha s_1$ (for $\alpha\in \CP^1-\{0,1,\infty\}$) define
smooth conics. 
Removing a neighborhood of a smooth fiber of the pencil generated by $s_0$
and $s_1$, we obtain a Lefschetz fibration over the disc, with fiber a sphere
with four punctures.
The three singular fibers of the pencil are nodal configurations consisting of
two tranversely intersecting spheres, with each component containing
two of the four base points; each of the three different
manners in which four points can be split into two groups of two is realized
at one of the singular fibers. The following diagram represents the three
singular conics of the pencil inside $\CP^2$ (left), and the corresponding
vanishing cycles inside a smooth fiber (right):
\medskip

\begin{center}
\epsfig{file=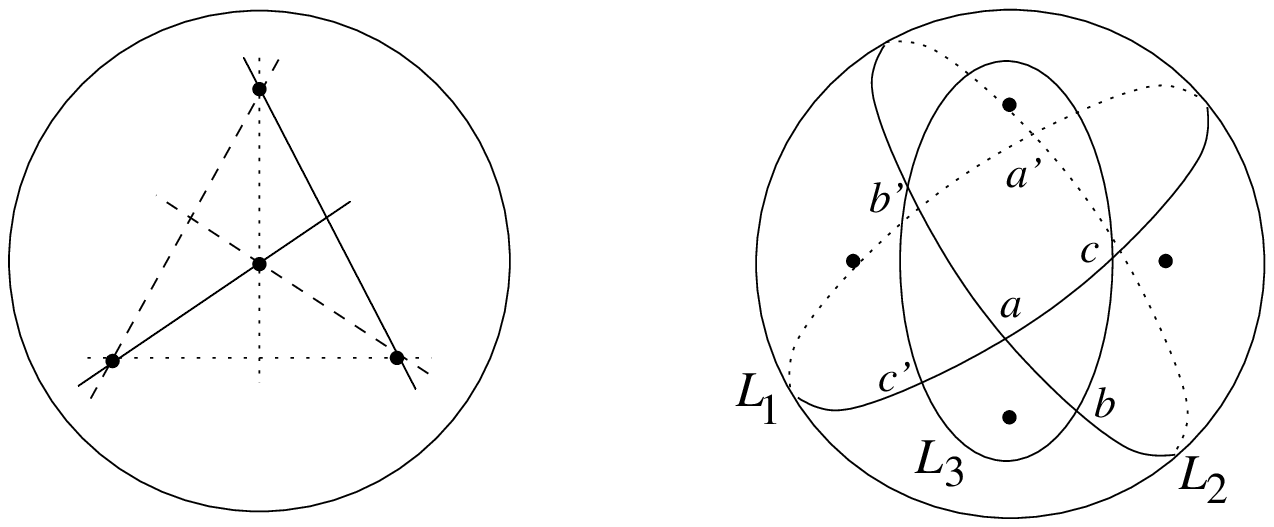,height=1.2in}
\end{center}

We can describe the monodromy of this Lefschetz pencil by
a homomorphism $\psi:\pi_1(\C-\{p_1,p_2,p_3\})\to \mathrm{Map}_{0,4}$ with
values in the mapping class group of a genus 0 surface with 4 boundary
components. After choosing
a suitable ordered basis of the free group $\pi_1(\C-\{p_1,p_2,p_3\})$,
we can make sure that $\psi$ maps the generators to the Dehn twists
$\tau_1,\tau_2,\tau_3$ along the three loops shown on the diagram.
On the other hand, because the normal bundles to the exceptional sections
of the blown-up Lefschetz fibration have degree $-1$, the monodromy at
infinity is given
by the boundary twist $\prod \delta_i$, the product of the four Dehn twists
along small loops encircling the four base points in the fiber; on the other
hand it is also the product of the monodromies
around each of the three singular fibers ($\tau_1,\tau_2,\tau_3$). Hence, the
monodromy of a pencil of conics in $\CP^2$ can be expressed by the relation
$\prod \delta_i=\tau_1\cdot \tau_2\cdot \tau_3$ in the mapping class group
$\mathrm{Map}_{0,4}$ ({\it lantern relation}).

Any two of the three vanishing cycles intersect transversely in two
points, so $\mathrm{Hom}(L_1,L_2)=\Z_2\,a\oplus \Z_2\,a'$,
$\mathrm{Hom}(L_2,L_3)=\Z_2\,b\oplus \Z_2\,b'$, and
$\mathrm{Hom}(L_1,L_3)=\Z_2\,c\oplus \Z_2\,c'$ are all two-dimensional.
There are no immersed 2-sided polygons in the punctured sphere $\Sigma_0$
with boundary in $L_i\cup L_j$ for any pair $(i,j)$, since each of the four
regions delimited by $L_i$ and $L_j$ contains one of the punctures, so
$\mu^1\equiv 0$. However, there are four triangles with boundary in
$L_1\cup L_2\cup L_3$ (with vertices $abc$, $ab'c'$, $a'b'c$, $a'bc'$
respectively), and in each case the cyclic ordering of the boundary is
compatible with the ordering of the vanishing cycles. Therefore, the
composition of morphisms is given by the formulas
$\mu^2(a,b)=\mu^2(a',b')=c$, $\mu^2(a,b')=\mu^2(a',b)=c'$.
Finally, the higher compositions $\mu^n$, $n\ge 3$ are all trivial in this
category, because the ordering condition $i_0<\dots<i_n$ never holds~\cite{Se2}.

\medskip

The objects $L_i$ of the category $\mathcal{F}_{vc}(f;\{\gamma_i\})$
actually correspond not only to Lagrangian spheres in
$\Sigma_0$ (the vanishing cycles), but also to Lagrangian discs in $X^0$
(the Lefschetz thimbles $D_i$); and the Floer intersection theory in
$\Sigma_0$ giving rise to $\mathrm{Hom}(L_i,L_j)$ and to the product
structures can also be thought of in terms of intersection theory for
the thimbles $D_i$ in $X^0$ (this is actually the reason of the
asymmetry between the cases $i<j$ and $i>j$ in Definition \ref{def:seidel}).
In any case, the properties of these objects depend very much
on the choice of the ordered collection of arcs $\{\gamma_i\}$.
Therefore, $\mathcal{F}_{vc}(f;\{\gamma_i\})$ has little geometric meaning
in itself, and should instead
be viewed as a collection of {\it generators} of a much larger
category which includes not only the Lefschetz thimbles, but also more general
Lagrangian submanifolds of $X^0$.
More precisely, the category naturally associated to the Lefschetz
pencil $f$ is not the finite directed $A_\infty$-category defined above, but rather the
(split-closed) {\it derived category} $D\mathcal{F}_{vc}(f)$ obtained from
$\mathcal{F}_{vc}(f;\{\gamma_i\})$ by considering (twisted) complexes
of formal direct sums of objects
(also including idempotent splittings and formal inverses
of quasi-isomorphisms). Replacing the ordered collection $\{\gamma_i\}$ by
another one $\{\gamma'_i\}$ leads to a different ``presentation'' of the
same derived category. Indeed, we have the following result \cite{Se1}:

\begin{theorem}[Seidel]
Given any two ordered collections $\{\gamma_i\}$ and
$\{\gamma'_i\}$, the categories $\mathcal{F}_{vc}(f;\{\gamma_i\})$ and
$\mathcal{F}_{vc}(f;\{\gamma'_i\})$ differ by a sequence of
{\it mutations} (operations that modify the ordering of the objects while
twisting some of them along others). Hence, the derived category
$D\mathcal{F}_{vc}(f)$ does not depend on the choice of $\{\gamma_i\}$.
\end{theorem}

Roughly speaking, complexes in the derived category correspond to
Lagrangian cycles obtained by gluing the corresponding Lefschetz thimbles.
For example, assume that the vanishing cycles $L_i$ and $L_j$ $(i<j)$ are
Hamiltonian isotopic to each other, so that (a smoothing of)
$\gamma_i\cup \gamma_j$ is a matching path. Then 
$\mathrm{Hom}(L_i,L_j)$ has rank two, with one generator in degree $0$ and
one in degree $n-1=\dim L_i$; let $a$ be the degree $0$ generator. Then
the complex $C=\{0\to L_i\stackrel{a}{\longrightarrow} L_j\to 0\}$, viewed as an
object of the derived category, represents the Lagrangian sphere associated
to the matching path $\gamma_i\cup \gamma_j$; for example it is easy to
check that $\mathrm{Hom}_{D\mathcal{F}_{vc}}(C,C)\simeq H^*(S^n,\Z_2)$.


Building Fukaya-type categories out of vanishing cycles may seem arbitrary,
but has a solid geometric underpinning, as suggested by the discussion
at the end of \S 5.1. In fairly general circumstances,
\emph{every} closed Lagrangian submanifold (with well-defined and
non-zero Floer homology with itself) of a K\"ahler manifold $X$ must intersect one of the vanishing
cycles of any Lefschetz pencil containing $X$ as a smooth member.  For the
theory of Biran and Cieliebak \cite{BC} shows that otherwise such a
submanifold could be displaced off itself by a Hamiltonian isotopy (first by
moving the Lagrangian into an open domain in $X$ admitting a
plurisubharmonic function with no top index critical points), a
contradiction to the non-triviality of Floer homology.  Seidel has pushed
this even further: for certain K3 surfaces arising as anticanonical divisors
in Fano 3-folds, every Lagrangian submanifold $L$ must have non-trivial
Floer cohomology with one of the vanishing cycles of the pencil.  Otherwise,
by repeatedly applying the exact sequence in Floer cohomology, one sees that
$HF(L,L)$ is graded isomorphic to a shifted version of itself, which is
absurd.  In this sense, the vanishing cycles really do "generate" the Fukaya
category; if $HF(L,V)=0$ for all vanishing cycles $V$, that is if $L$ has no
morphisms to any vanishing cycle, then $HF(L,L)=0$ so $L$ has no identity
morphism, and must represent the zero object of the category.


\subsection{Applications to mirror symmetry}

The construction described above has various applications to homological
mirror symmetry. In the context of Calabi-Yau manifolds, these have
to do with a conjecture of Seidel about the relationship between the derived
Fukaya category $D\mathcal{F}_{vc}(f)$ of the Lefschetz pencil $f$ and the derived
Fukaya category $D\mathcal{F}(X)$ of the closed symplectic manifold $X$
\cite{SeICM}. As seen above, when passing to the derived category
$D\mathcal{F}_{vc}(f)$, we hugely increase the number of objects, by considering not
only the thimbles $D_i$ but also arbitrary complexes obtained from them;
this means that the objects of $D\mathcal{F}_{vc}(f)$ include all sorts of
(not necessarily closed) Lagrangian submanifolds in $X^0$, with boundary
in $\Sigma_0$. Since Fukaya categories are only concerned with closed
Lagrangian submanifolds, it is necessary to consider a subcategory of
$D\mathcal{F}_{vc}(f)$ whose objects correspond to the closed Lagrangian
submanifolds in $X^0$ (i.e., combinations of $D_i$ for which the
boundaries cancel); it is expected that this can be done in purely
categorical terms by considering those objects of $D\mathcal{F}_{vc}(f)$ on which
the Serre functor acts simply by a shift. The resulting subcategory
should be closely related to the derived Fukaya category of the open
manifold $X^0$. This leaves us with the problem of relating
$\mathcal{F}(X^0)$ with $\mathcal{F}(X)$. These two categories have
the same objects and morphisms (Lagrangians in $X$ can be made disjoint
from $\Sigma_\infty$), but the differentials and product structures differ.
More precisely, the definition of $\mu^n$ in $\mathcal{F}(X^0)$ only
involves counting pseudo-holomorphic discs contained in $X^0$, i.e.\ disjoint
from the hypersurface $\Sigma_\infty$. In order to account for the missing
contributions, one should introduce a formal parameter $q$ and count the
pseudo-holomorphic discs with boundary in $\bigcup L_i$ that intersect
$\Sigma_\infty$ in $m$ points (with multiplicities) with a coefficient
$q^m$. The introduction of this parameter $q$ leads to a {\it deformation of
$A_\infty$-structures}, i.e.\ an
$A_\infty$-category in which the differentials and products $\mu^n$ are
defined over a ring of
formal power series in the variable $q$; the limit $q=0$ corresponds to
the (derived) Fukaya category $D\mathcal{F}(X^0)$, while non-zero values of
$q$ are expected to yield $D\mathcal{F}(X)$.

These considerations provide a strategy to calculate Fukaya categories
(at least for some examples) by induction on dimension \cite{SeICM};
an important recent development in this direction is Seidel's proof of
homological mirror symmetry for quartic K3 surfaces \cite{Se4}.
\medskip

Returning to more elementary considerations,
another context in which the construction of Definition \ref{def:seidel}
is relevant is that of mirror
symmetry for Fano manifolds ($c_1(TM)>0$). Rather than manifolds,
mirrors of Fano manifolds are {\it Landau-Ginzburg models}, i.e.\ pairs
$(Y,w)$, where $Y$ is a non-compact manifold and $w:Y\to\C$, the
``superpotential'', is a holomorphic function. The complex (resp.\
symplectic) geometry of a Fano manifold $M$ is then expected to correspond
to the symplectic (resp.\ complex) geometry of the {\it critical points} of
the superpotential $w$ on its mirror. In particular, the homological mirror
conjecture of Kontsevich now asserts that the categories $D^b\mathrm{Coh}(M)$
and $D\mathcal{F}_{vc}(w)$ should be equivalent. 
Following the ideas of various
people (Kontsevich, Hori, Vafa, Seidel, etc.), the conjecture can be
verified on many examples, at least in those cases where the critical
points of $w$ are isolated and non-degenerate. 

For example, let $a,b,c$ be three mutually prime positive integers, and
consider the weighted projective plane $M=\CP^2(a,b,c)=(\C^3-\{0\})/\C^*$,
where $\C^*$ acts by $t\cdot (x,y,z)=(t^a x,t^b y,t^c z)$. In general $M$
is a Fano orbifold (when $a=b=c=1$ it is the usual projective plane).
Consider the mirror Landau-Ginzburg model $(Y,w)$, where $Y=\{x^ay^bz^c=1\}$
is a hypersurface in $(\C^*)^3$, equipped with an exact K\"ahler form, and
$w=x+y+z$. The superpotential $w$ has $a+b+c$ isolated non-degenerate
critical points, and hence determines an affine Lefschetz fibration to which
we can apply the construction of Definition \ref{def:seidel}.
To be more precise, in this context one should actually
define the Fukaya category of vanishing cycles over a coefficient
ring $R$ larger than $\Z_2$ or $\Z$, for example $R=\C$, counting each
pseudo-holomorphic disc $u:D^2\to \Sigma_0$ with a weight
$\pm \exp(-\int_{D^2} u^*\omega)$, or a Novikov ring.
Then we have the following result:

\begin{theorem}[\cite{AKO}]
The categories $D^b\mathrm{Coh}(M)$ and $D\mathcal{F}_{vc}(w)$ are
equivalent.
\end{theorem}

The proof relies on the identification of suitable collections of generators
for both categories, and an explicit verification that the morphisms on both
sides obey the same composition rules. Moreover, it can also be shown
that non-exact deformations of the K\"ahler structure on $Y$ correspond
to non-commutative deformations of $M$, with an explicit relationship
between the deformation parameter on $M$ and the complexified K\"ahler
class $[\omega+iB]\in H^2(Y,\C)$ \cite{AKO}.

The homological mirror conjecture can be similarly verified for
various other examples, and one may reasonably expect that, in
the near future, vanishing cycles will play
an important role in our understanding of mirror symmetry, not only for
Fano and Calabi-Yau varieties, but maybe also for varieties of general~type.


\begin{thebibliography}{ADKY}
\bibitem[ABKP]{ABKP}
J. Amor\'os, F. Bogomolov, L. Katzarkov, T. Pantev, 
{\sl Symplectic Lefschetz fibrations with arbitrary fundamental groups}, 
J. Differential Geom. {\bf 54} (2000), 489--545.
\bibitem[AGTV]{AGTV}
M. Amram, D. Goldberg, M. Teicher, U. Vishne, {\sl The fundamental group
of the Galois cover of the surface $\CP^1\times T$}, Alg. Geom. Topol.
{\bf 2} (2002), 403--432.
\bibitem[Au1]{Au1}
D.\ Auroux, {\sl Asymptotically holomorphic families of
symplectic submanifolds}, Geom.\ Funct.\ Anal. {\bf 7} (1997), 971--995.
\bibitem[Au2]{Au2} D. Auroux,
{\sl Symplectic 4-manifolds as branched coverings of $\mathbb{CP}^2$},
Invent. Math. {\bf 139} (2000), 551--602.
\bibitem[Au3]{Au3} D. Auroux, 
{\sl Symplectic maps to projective spaces and symplectic invariants},
Turkish J.\ Math.\ {\bf 25} (2001), 1--42 ({\tt math.GT/0007130}).
\bibitem[Au4]{Au4} D. Auroux,
{\sl Estimated transversality in symplectic geometry and projective maps},
``Symplectic Geometry and Mirror Symmetry'', Proc. 4th KIAS International
Conference, Seoul (2000), World Sci., Singapore, 2001, pp. 1--30
({\tt math.SG/0010052}).
\bibitem[Au5]{Au5}
D. Auroux, {\sl A remark about Donaldson's construction of symplectic
submanifolds}, J.\ Symplectic Geom.\ {\bf 1} (2002), 647--658
({\tt math.DG/0204286}).
\bibitem[Au6]{Au6}
D. Auroux, {\sl Fiber sums of genus 2 Lefschetz fibrations},
Turkish J.\ Math.\ {\bf 27} (2003), 1--10 ({\tt math.GT/0204285}).
\bibitem[ADK]{ADK}
D. Auroux, S.\,K. Donaldson, L. Katzarkov,
{\sl Luttinger surgery along Lagrangian tori and non-isotopy for singular 
symplectic plane curves}, Math. Ann. {\bf 326} (2003), 185--203.
\bibitem[ADKY]{ADKY}
D. Auroux, S.\,K. Donaldson, L. Katzarkov, M. Yotov,
{\sl Fundamental groups of complements of plane curves and symplectic 
invariants}, preprint, to appear in Topology ({\tt math.GT/0203183}).
\bibitem[AK]{AK}
D. Auroux, L. Katzarkov,
{\sl Branched coverings of $\mathbb{CP}^2$ and invariants of symplectic 
4-manifolds}, Invent.\ Math. {\bf 142} (2000), 631--673.
\bibitem[AK2]{AK2}
D. Auroux, L. Katzarkov, {\sl The degree doubling formula
for braid monodromies and Lefschetz pencils}, preprint.
\bibitem[AKO]{AKO}
D. Auroux, L. Katzarkov, D. Orlov, {\sl Homological mirror symmetry for
weighted projective planes}, in preparation.
\bibitem[AMP]{AMP}
D. Auroux, V. Mu\~noz, F. Presas, {\sl Lagrangian submanifolds and Lefschetz
pencils}, in preparation.
\bibitem[Bi]{Biran}
P. Biran, {\sl A stability property of symplectic packing in dimension 4},
Invent. Math. {\bf 136} (1999), 123--155.
\bibitem[BC]{BC}
P. Biran, K. Cieliebak, {\sl Symplectic topology on subcritical manifolds}, 
Comment. Math. Helv. {\bf 76} (2001), 712-753.
\bibitem[Bir]{Bi}
J. Birman, {\it Braids, Links and Mapping class groups}, Annals of
Math. Studies {\bf 82}, Princeton Univ. Press, Princeton, 1974.
\bibitem[CMT]{CMT}C. Ciliberto, R. Miranda, M. Teicher, {\sl Pillow
degenerations of K3 surfaces}, Applications of Algebraic Geometry to Coding
Theory, Physics, and Computation, NATO Science Series II, vol.\ {\bf 36},
2001, pp.\ 53--63.
\bibitem[Do1]{Do1} S.\,K.\ Donaldson, {\sl Symplectic submanifolds and
almost-complex geometry}, J.\ Differential Geom. {\bf 44} (1996), 666--705.
\bibitem[Do2]{Do2} S.K. Donaldson, {\sl Lefschetz fibrations in symplectic
geometry}, Documenta Math., Extra Volume ICM 1998, II, 309--314.
\bibitem[Do3]{Do3} S.K. Donaldson, {\sl Lefschetz pencils on symplectic
manifolds}, J. Differential Geom. {\bf 53} (1999), 205--236.
\bibitem[DS]{DS}
S. Donaldson, I. Smith, {\sl Lefschetz pencils and the
canonical class for symplectic 4-manifolds}, Topology {\bf 42} (2003),
743--785.
\bibitem[FO${}^3$]{FO3} K.\ Fukaya, Y.-G.\ Oh, H.\ Ohta, K.\ Ono, {\sl
Lagrangian intersection Floer theory: Anomaly and obstruction}, preprint.
\bibitem[FS]{FS2}
R.\ Fintushel, R.\ Stern, {\sl Symplectic surfaces in a fixed homology
class}, J.\ Differential Geom.\ {\bf 52} (1999), 203--222.
\bibitem[Gi]{Gi}
E. Giroux, {\sl G\'eom\'etrie de contact: de la
dimension trois vers les dimensions sup\'erieures},
Proc.\ International Congress of Mathematicians, Vol. II (Beijing, 2002),
Higher Ed.\ Press, Beijing, 2002, pp. 405--414.
\bibitem[Go1]{Go1} R.\,E. Gompf, {\sl A new construction of symplectic
manifolds}, Ann. Math. {\bf 142} (1995), 527--595.
\bibitem[Go2]{Go2} R.\,E. Gompf, {\sl A topological characterization of
symplectic manifolds}, preprint ({\tt math.SG/ 0210103}).
\bibitem[GS]{GS} R.\,E. Gompf, A.\,I. Stipsicz, {\it 4-manifolds and Kirby
calculus}, Graduate Studies in Math. {\bf 20}, Amer. Math. Soc., Providence,
1999.
\bibitem[Gr]{Gr} M.\ Gromov, {\sl Pseudo-holomorphic curves in symplectic
manifolds}, Invent.\ Math. {\bf 82} (1985), 307--347.
\bibitem[IMP]{IMP} A. Ibort, D. Martinez-Torres, F. Presas,
{\sl On the construction of contact submanifolds with
prescribed topology}, J.\ Differential Geom. {\bf 56} (2000), 235--283.
\bibitem[KK]{KK} V. Kharlamov, V. Kulikov, {\sl On braid monodromy
factorizations}, Izvestia Math.\ {\bf 67} (2003), 79--118
({\tt math.AG/0302113}).
\bibitem[Ku]{Ku} V. Kulikov, {\sl On a Chisini conjecture}, Izvestia Math.\
{\bf 63} (1999), 1139--1170. 
\bibitem[McS]{McS2} D.\ McDuff and D.\ Salamon, {\it J-holomorphic curves
and quantum cohomology}, Univ.\ Lecture Series No.\ {\bf 6},
Amer.\ Math.\ Soc., Providence, 1994.
\bibitem[Mo1]{Mo1}
B. Moishezon, {\sl Complex surfaces and connected sums of complex projective
planes}, Lecture Notes in Math.\ {\bf 603}, Springer, Heidelberg, 1977.
\bibitem[Mo2]{Mo2} 
B. Moishezon, {\sl Stable branch curves and braid
monodromies}, Algebraic Geometry (Chicago, 1980), Lecture Notes in Math. 
{\bf 862}, Springer, Heidelberg, 1981, pp.\ 107--192.
\bibitem[Mo3]{Mo3}
B.\ Moishezon, {\sl On cuspidal branch curves}, J. Algebraic Geom. {\bf 2}
(1993), 309--384.
\bibitem[Mo4]{Mo4}
B.\ Moishezon, {\sl The arithmetic of braids and a
statement of Chisini}, Geometric Topology (Haifa, 1992), Contemp. Math.
{\bf 164}, Amer. Math. Soc., Providence, 1994, pp.\ 151--175.
\bibitem[Mo5]{MVeronese}
B.\ Moishezon, {\sl Topology of generic polynomial maps in complex dimension
two}, preprint.
\bibitem[MRT]{MRT}
B.\ Moishezon, A.\ Robb, M.\ Teicher, {\sl On Galois covers of Hirzebruch
surfaces}, Math. Ann. {\bf 305} (1996), 493--539.
\bibitem[Na]{Nakajima}
H.\ Nakajima, {\sl Lectures on Hilbert schemes of points on surfaces}, Univ.\ Lecture Series No.\ {\bf 18},
Amer.\ Math.\ Soc., Providence, 1999.
\bibitem[Oz]{Oz}
B. Ozbagci, {\sl Signatures of Lefschetz fibrations}, Pacific J.\ Math.
{\bf 202} (2002), 99--118.
\bibitem[OS]{OS}
B. Ozbagci, A. Stipsicz, {\sl Noncomplex smooth 4-manifolds with
genus 2 Lefschetz fibrations}, Proc.\ Amer.\ Math.\ Soc.\ {\bf 128} (2000),
3125--3128.
\bibitem[Ro]{Ro}
A.\ Robb, {\sl On branch curves of algebraic surfaces}, Singularities and
Complex Geometry (Beijing, 1994), Amer.\ Math.\ Soc./Int.\ Press\,\ %
Stud.\ Adv.\ Math. {\bf 5}, Amer.\ Math.\ Soc., Providence, 1997, pp.\
193--221.
\bibitem[Se1]{Se1}
P. Seidel, {\sl Vanishing cycles and mutation},
Proc. 3rd European Congress of Mathe\-ma\-tics (Barcelona, 2000), Vol. II,
Progr.\ Math.\ {\bf 202}, Birkh\"auser, Basel, 2001, pp.\ 65--85
({\tt math.SG/0007115}).
\bibitem[Se2]{Se2}
P. Seidel, {\sl More about vanishing cycles and mutation},
``Symplectic Geometry and Mirror Symmetry'', Proc. 4th KIAS International
Conference, Seoul (2000), World Sci., Singapore, 2001, pp. 429--465
({\tt math.SG/0010032}).
\bibitem[Se3]{SeICM}
P. Seidel, {\sl Fukaya categories and deformations},
Proc.\ International Congress of Mathematicians, Vol. II (Beijing, 2002),
Higher Ed.\ Press, Beijing, 2002, pp. 351--360.
\bibitem[Se4]{Se3}
P. Seidel, {\sl Fukaya categories and Picard-Lefschetz theory}, in
preparation.
\bibitem[Se5]{Se4}
P. Seidel, {\sl Homological mirror symmetry for the quartic surface},
preprint ({\tt math.SG/} {\tt 0310414}).
\bibitem[Sh]{Sh}
V. Shevchishin, {\sl On the local version of the Severi problem}, preprint
({\tt math.AG/0207048}).
\bibitem[ST]{ST}
B. Siebert, G. Tian,
{\sl On the holomorphicity of genus two Lefschetz fibrations}, preprint,
to appear in Ann.\ Math ({\tt math.SG/0305343}).
\bibitem[Sm1]{SmHodge}
I. Smith, {\sl Lefschetz fibrations and the Hodge bundle},
Geom.\ Topol. {\bf 3} (1999), 211--233.
\bibitem[Sm2]{Sm1}
I. Smith, {\sl Lefschetz pencils and divisors in moduli
space}, Geom.\ Topol. {\bf 5} (2001), 579--608. 
\bibitem[Sm3]{Sm2}
I. Smith, {\sl Geometric monodromy and the hyperbolic disc},
Quarterly J.\ Math.\ {\bf 52} (2001), 217--228 ({\tt math.SG/0011223}).
\bibitem[Sm4]{Sm3}
I. Smith, {\sl Serre-Taubes duality for pseudoholomorphic curves},
Topology {\bf 42} (2003), 931--979.
\bibitem[Ta1]{Ta1} C.\,H.\ Taubes, {\sl The Seiberg-Witten and the Gromov
invariants}, Math.\ Res.\ Lett.\ {\bf 2} (1995), 221--238.
\bibitem[Ta2]{Ta2} C.\,H.\ Taubes, {\sl The geometry of the Seiberg-Witten
invariants}, Surveys in Differential Geometry, Vol.\ III (Cambridge, 1996),
Int.\ Press, Boston, 1998, pp.\ 299--339.
\bibitem[Te1]{Te1}
M. Teicher, {\sl Braid groups, algebraic surfaces and
fundamental groups of complements of branch curves}, Algebraic Geometry
(Santa Cruz, 1995), Proc. Sympos. Pure Math., {\bf 62} (part 1), Amer. 
Math. Soc., Providence, 1997, pp.\ 127--150.
\bibitem[Te2]{Te2}
M.\ Teicher, {\sl New invariants for surfaces}, Tel Aviv Topology Conference:
Rothenberg Festschrift (1998), Contemp. Math. {\bf 231}, Amer. Math. Soc.,
Providence, 1999, pp.\ 271--281 ({\tt math.AG/9902152}).
\bibitem[Te3]{TVeronese}
M.\ Teicher, {\sl The fundamental group of a $\CP^2$ complement of
a branch curve as an extension of a solvable group by a symmetric group},
Math. Ann. {\bf 314} (1999) 19--38.
\bibitem[Th]{Th} W.\ Thurston, {\sl Some simple examples of symplectic
manifolds}, Proc.\ Amer.\ Math.\ Soc. {\bf 55} (1976), 467--468.
\bibitem[Us]{Usher}
M.\ Usher, {\sl The Gromov invariant and the Donaldson-Smith standard surface count},
preprint ({\tt math.SG/0310450}).
\end{thebibliography}
\end{document}